\documentclass{article}
\title{On Higher Frobenius-Schur Indicators}
\author{Yevgenia Kashina \quad Yorck Sommerh\"auser \vspace{3mm}\\  
Yongchang Zhu}
\date{}

\usepackage{amsmath}
\usepackage{theorem}
\usepackage{amssymb}
\usepackage{makeidx}
\usepackage{latexsym}
\usepackage{diagramb}
\usepackage[all]{xy}

\SilentMatrices

\makeindex

\makeatletter
\renewcommand{\subsection}{\@startsection{subsection}{2}{0em}%
{\baselineskip}{-0em}{\bfseries\normalsize}}
\makeatother

\makeatletter
\newcommand{\listofdefinitions}{\@starttoc{def}}

\newcommand{\l@definition}[2]{\par\noindent#1 {\itshape #2}}
\makeatother

\theoremstyle{plain}
\theorembodyfont{\rmfamily}

\newtheorem{thm}{Theorem}

\newtheorem{prop}[thm]{Proposition}
\newtheorem{lemma}[thm]{Lemma}
\newtheorem{corollary}[thm]{Corollary}

\newtheorem{pf}{Proof.}

\newtheorem{defn}[thm]{Definition}

\theoremstyle{break}
\theorembodyfont{\rmfamily}

\newcommand{\qed}{$\Box$}

\newcommand{\rank}{\operatorname{rank}}

\newcommand{\Ch}{\operatorname{Ch}}
\newcommand{\Gal}{\operatorname{Gal}}
\newcommand{\End}{\operatorname{End}}

\newcommand{\Span}{\operatorname{Span}}

\newcommand{\Hom}{\operatorname{Hom}}
\newcommand{\id}{\operatorname{id}}

\newcommand{\Tr}{\operatorname{tr}}

\newcommand{\diag}{\operatorname{diag}}
\newcommand{\op}{\scriptstyle \operatorname{op}}
\newcommand{\cop}{\scriptstyle \operatorname{cop}}

\newcommand{\ind}{\operatorname{ind}}
\newcommand{\ord}{\operatorname{ord}}
\newcommand{\mult}{\operatorname{mult}}
\newcommand{\Res}{\operatorname{Res}}

\def\1{{(1)}}
\def\2{{(2)}}
\def\3{{(3)}}
\def\4{{(4)}}
\def\5{{(5)}}
\def\6{{(6)}}
\def\7{{(7)}}
\def\8{{(8)}}
\def\9{{(9)}}

\def\m{{(m)}}

\def\o{\otimes}
\def\l{\langle}
\def\r{\rangle}

\def\ra{\rightarrow}
\def\N{{\mathbb N}}
\def\Z{{\mathbb Z}}
\def\Q{{\mathbb Q}}
\def\R{\mathbb R}
\def\C{\mathbb C}

\begin{document}

\maketitle

\begin{abstract}
\hspace{-5mm}We study the higher Frobenius-Schur indicators of modules over
semisimple Hopf algebras, and relate them to other invariants as the exponent, the order, and the index. We prove various divisibility and integrality results for these invariants. In particular, we prove a version of Cauchy's theorem for semisimple Hopf algebras. Furthermore, we give some examples that illustrate the general theory.
\end{abstract}


\section*{Introduction} \label{Sec:Introd}
\addcontentsline{toc}{section}{Introduction}
For a finite group, one can evaluate a character on the sum of all $m$-th powers of the group elements. The resulting number, divided by the order of the group, is called the $m$-th Frobenius-Schur indicator of the character. The first use of these indicators was made by F.~G.~Frobenius and I.~Schur (cf.~\cite{FS}) to give a criterion when a representation of a finite group can be realized by matrices with real entries---for this question, it is the second indicator that is relevant. This is also meaningful for other fields than the complex numbers: Here the indicator tells whether or not a given module is self-dual. 

Higher indicators, i.e., indicators with $m>2$, arise when one considers
the root number function in a finite group. This function assigns to a group element the number of its $m$-th roots, i.e., the number of group elements whose $m$-th power is equal to the given element. It is clear that this number depends only on the conjugacy class, and therefore defines a class function that can be expanded in terms of the irreducible characters. Using the orthogonality relations for characters, it is not hard to see that the coefficient of an irreducible character in this expansion is its $m$-th Frobenius-Schur indicator (cf.~\cite{Isaacs}, Lem.~(4.4), p.~49).

For Hopf algebras, Frobenius-Schur indicators were first considered by V.~Lin\-chen\-ko and S.~Montgomery on the one hand (cf.~\cite{LinMont1}) and by J.~Fuchs, A.~Ch. Gan\-chev, K.~Szlach\'{a}nyi, and P.~Vecserny\'{e}s on the other hand (cf.~\cite{FuGanSzlVec1}). Here, the sum of the $m$-th powers of the group elements is replaced by the $m$-th Sweedler power of the integral. The authors then use the indicators, or at least the second indicator, to prove an analogue of the criterion of Frobenius and Schur whether or not a representation is self-dual: The Frobenius-Schur theorem asserts that this depends on whether the second indicator is $0$, $1$, or $-1$.

The topic of the present writing are the higher Frobenius-Schur indicators for semisimple Hopf algebras and their relation to other invariants of irreducible characters. These other invariants are the order, the multiplicity, the exponent, and the index. Let us briefly describe the nature of these invariants. The notion of the order of an irreducible character is a generalization of the notion of the order of an element in a finite group: It is the smallest  integer such that the corresponding tensor power contains a nonzero invariant subspace. The dimension of this invariant subspace is called the multiplicity of the irreducible character. An irreducible character has order~1 if and only if it is trivial, and has order~2 if and only if it is self-dual. In these cases, the multiplicity of the character is~1.

The exponent of a semisimple Hopf algebra is another invariant, which generalizes the exponent of a group (cf.~\cite{KashPow}). The exponent of a module is a slight generalization
of this concept: In the group case, it is the exponent of the image of the group in the representation. There are various ways to generalize this concept to semisimple Hopf algebras; we will use Sweedler powers on the one hand and a certain canonical tensor on the other hand.

The next invariant that we study, the index of imprimitivity, arises from Perron-Frobenius theory. It is clear that the matrix representation of the left multiplication by the character of a module with respect to the basis consisting of all irreducible characters has nonnegative integer entries. As we will explain below, the corresponding Perron-Frobenius eigenvalue, i.e., the positive eigenvalue that has the largest possible absolute value, is the degree of the character. However, since the entries of the above matrix are in general not strictly positive, this eigenvalue is not necessarily strictly greater than the absolute values of the other characters, so that there can be other eigenvalues which are not positive, but have the same absolute value. As we will see, in the most interesting cases the above matrix is indecomposable; in this case, the number of such eigenvalues is called the index of imprimitivity.  

The text is organized as follows: In Section~\ref{Sec:SweedlPow}, we discuss the formalism of Sweedler powers in an arbitrary bialgebra. A Sweedler power of an element in a bialgebra is constructed by applying the comultiplication several times, permuting the arising tensor factors and multiplying them together afterwards. This notion is a slight modification of the original notion (cf.~\cite{KashPow}), where the tensor factors were not permuted, and has the advantage that iterated Sweedler powers are still Sweedler powers. We then consider the Sweedler powers that arise from a certain special kind of permutations. This is motivated by the fact that the values of characters on the Sweedler powers of the integral lie in certain cyclotomic fields, and this special kind of Sweedler powers is well adapted to describe the action of the Galois group on these values.

From Section~\ref{Sec:FrobSchur} on, we consider semisimple Hopf algebras over algebraically closed fields of characteristic zero. 
We prove a first formula for the higher Frobenius-Schur indicators that should be understood as a generalization of the Frobenius-Schur theorem for these indicators---in particular, it implies the Frobenius-Schur theorem for the second indicators. This first formula describes Frobenius-Schur indicators in terms of a certain operator on the corresponding tensor power of the module, and we establish several other properties of this operator as well. 

In Section~\ref{Sec:Exponent}, we then consider the exponent of a module and prove a second formula for the higher Frobenius-Schur indicators that uses a certain canonical tensor. Combining this with the first formula, we prove a version of Cauchy's theorem for Hopf algebras: A prime that divides the dimension of a semisimple Hopf algebra must also divide its exponent. This result was conjectured by P.~Etingof and S.~Gelaki (cf.~\cite{EG1}); it was known in the case of the prime~2 (cf.~\cite{KSZ}). Furthermore, we prove that the higher indicators are integers if the exponent is squarefree.

In Section~\ref{Sec:Order}, we define the notion of the order and the multiplicity of a module and prove that the order of a module divides its multiplicity times the dimension of the Hopf algebra. This result generalizes the theorem that a semisimple Hopf algebra that has a nontrivial self-dual simple module must have even dimension (cf.~\cite{KSZ}) to modules of arbitrary orders---nontrivial self-dual simple modules are of order~2, and the multiplicity is~1 in this case. In particular, we get in this way a fully independent new proof of the old theorem.

In Section~\ref{Sec:Index}, we study the index of imprimitivity, or briefly the index, of the matrix that represents the left multiplication by a character with respect to the canonical basis that we have in the character ring---the basis consisting of the irreducible characters. The main result of this section is a precise formula for the index in terms of central grouplike elements. Essentially, the result says that the eigenvalues
of the above matrix that have the same absolute value as the degree
are obtained by evaluating the character at certain central grouplike elements. As a consequence of this formula, we see that the index divides the order as well as the exponent. 

In Section~\ref{Sec:DrinfDouble}, we apply a new tool---the Drinfel'd double of the Hopf algebra. We prove that, by restricting modules over the
Drinfel'd double to the Hopf algebra, we get a map from the character ring
of the Drinfel'd double onto the center of the character ring of the Hopf algebra. From this, we deduce that the center of the rational character  ring of the Hopf algebra, i.e., the span of the irreducible characters over the rational numbers, is isomorphic to a product of subfields of the cyclotomic field whose order is the exponent of the Hopf algebra. Finally,
we deduce a third formula for the Frobenius-Schur indicators in terms of the action of the Drinfel'd element on the induced module over the Drinfel'd double.

In Section~\ref{Sec:Examples}, we finally compute explicitly a number of examples. In this way, we can limit the possible generalizations of the results that we have obtained. The class of examples that we study are certain extensions of group rings by dual group rings. In particular, the Drinfel'd doubles of finite groups belong to this class.

Throughout the whole exposition, we consider a base field that is denoted by~$K$. All vector spaces considered are defined over~$K$, and all tensor products without subscripts are taken over~$K$. Unless stated otherwise, a module is a left module. The set of natural numbers is the set~$\N:=\{1,2,3,\ldots\}$; in particular, $0$ is not a natural number. 
The symbol~$\Q_n$\index{$\Q_n$} denotes the $n$-th cyclotomic field, and not the field of $n$-adic numbers, and~$\Z_n$\index{$\Z_n$} denotes the set~$\Z/n\Z$ of integers modulo~$n$, and not the ring of $n$-adic integers.
 
From Section~\ref{Sec:FrobSchur} on, except for Corollary~\ref{Cauchy}, we assume that the base field~$K$ is algebraically closed of characteristic zero. $H$ denotes a semisimple Hopf algebra with coproduct~$\Delta$\index{$\Delta$}, counit~$\varepsilon$\index{$\varepsilon$}, and antipode~$S$\index{$S$}. We will use the same symbols to denote the corresponding structure elements of the dual Hopf algebra~$H^*$. Note that a semisimple Hopf algebra is automatically finite-dimensional (cf.~\cite{Sw}, Chap.~V, Exerc.~4, p.~108). By results of R.~G.~Larson and D.~E.~Radford (cf.~\cite{LR2}, Thm.~3.3, p.~276; \cite{LR1}, Thm.~4, p.~195), $H$ is also cosemisimple and involutory, i.e., $H^*$ is semisimple and the antipode is an involution. From Maschke's theorem (cf.~\cite{M}, Thm.~2.2.1, p.~20), we get that there is a unique two-sided integral $\Lambda \in H$ such that $\varepsilon(\Lambda)=1$; this\index{$\Lambda$} element will be used heavily throughout.

Furthermore, we use the convention that propositions, definitions, and similar items are referenced by the paragraph in which they occur; they are only numbered separately if this reference is ambiguous.

\section{The Calculus of Sweedler Powers} \label{Sec:SweedlPow}
\subsection[Monotone maps]{} \label{Monoton}
Let us begin by introducing some notation that will be used throughout this section. For a natural number~$n$, which is by definition at least~1, we use the notation
$I_n:= \{1,\ldots,n\}$\index{$I_n$}. Suppose now that $m_1,m_2,\ldots,m_k$ are several natural numbers. On the product
$I_{m_1} \times I_{m_2} \times \ldots \times I_{m_k}$,
we introduce the so-called lexicographical ordering. This means that, for two $k$-tuples $(i_1,i_2,\ldots,i_k)$ and $(j_1,j_2,\ldots,j_k)$ in this set, we define $(i_1,i_2,\ldots,i_k) < (j_1,j_2,\ldots,j_k)$ if and only if there is an index $l \leq k$ such that $i_1=j_1, i_2=j_2, \ldots, i_{l-1}=j_{l-1}$, but $i_l<j_l$. As this is a total ordering, there is a unique strictly monotone map
$$\varphi_{m_1,\ldots,m_k}: I_{m_1} \times I_{m_2} \times \ldots \times I_{m_k} \longrightarrow I_n$$
\index{$\varphi_{m_1,\ldots,m_k}$}
where $n := m_1 \cdot m_2 \cdot \ldots \cdot m_k$. Explicitly, this map is given by the formula
$$\varphi_{m_1,\ldots,m_k}(i_1,\ldots,i_k) = 
(i_1-1)n_2 + (i_2-1) n_3 + \ldots + (i_{k-1}-1) n_k + i_k$$
where $n_i:=m_i \cdot m_{i+1} \cdot \ldots \cdot m_k$, so that $n_k=m_k$ and $n_1=n$.
This holds since if
$i_r=j_r$ for $r=1,\ldots,l-1$, but $i_l<j_l$, then we have 
$$(i_l-1) n_{l+1} + \ldots + (i_{k-1}-1) n_k + i_k
< (j_l-1) n_{l+1} + \ldots + (j_{k-1}-1) n_k + j_k$$
as the maximal value that $(i_{l+1}-1) n_{l+2}  + \ldots + (i_{k-1}-1) n_k + i_k$ can attain is
$$(m_{l+1}-1) n_{l+2}  + \ldots + (m_{k-1}-1) n_k + m_k = m_{l+1} n_{l+2}  = n_{l+1}$$
and $j_k \geq 1$. Considering this calculation for $l=0$, we also see that the map is well-defined. Obviously, we have $\varphi_n = \id_{I_n}$ if $k=1$.

Now suppose that, for a different~$l$, we have another set of natural numbers $m'_1,m'_2,\ldots,m'_l$, and let 
$n':=m'_1 \cdot m'_2 \cdot \ldots \cdot m'_l$. On the one hand, we can then look at the composition
$$(I_{m_1} \times \ldots \times I_{m_k}) \times 
(I_{m'_1} \times \ldots \times I_{m'_l})
\xrightarrow{\varphi_{m_1,\ldots,m_k} \times \varphi_{m'_1,\ldots,m'_l}} I_n \times I_{n'} \xrightarrow{\varphi_{n,n'}} I_{nn'}$$
On the other hand, we have the map
$$\varphi_{m_1,\ldots,m_k, m'_1,\ldots,m'_l}:
I_{m_1} \times \ldots \times I_{m_k} \times 
I_{m'_1} \times \ldots \times I_{m'_l} \longrightarrow I_{nn'}$$
It is clear from the definition of the lexicographical ordering that we have
$$(i_1,\ldots,i_k,i'_1,\ldots,i'_l) < (j_1,\ldots,j_k,j'_1,\ldots,j'_l)$$
if and only if we have $(i_1,\ldots,i_k) < (j_1,\ldots,j_k)$
or $(i_1,\ldots,i_k) = (j_1,\ldots,j_k)$ and $(i'_1,\ldots,i'_l) < (j'_1,\ldots,j'_l)$. This shows that, if all appearing product sets are ordered lexicographically, both maps considered above are strictly monotone, and therefore must be equal. This proves the following statement:
\begin{lemma}
For natural numbers $m_1,\ldots,m_k$ and $m'_1,\ldots,m'_l$, we have
$$\varphi_{m_1,\ldots,m_k, m'_1,\ldots,m'_l} =
\varphi_{n,n'} \circ (\varphi_{m_1,\ldots,m_k} \times \varphi_{m'_1,\ldots,m'_l})$$
where $n := m_1 \cdot m_2 \cdot \ldots \cdot m_k$ and 
$n':=m'_1 \cdot m'_2 \cdot \ldots \cdot m'_l$.
\end{lemma}
Of course, this can also be seen from the explicit formulas, as we have
\begin{align*}
\varphi_{n,n'}(&\varphi_{m_1,\ldots,m_k}(i_1,\ldots,i_k),
\varphi_{m'_1,\ldots,m'_l}(i'_1,\ldots,i'_l))\\
&= n'(\varphi_{m_1,\ldots,m_k}(i_1,\ldots,i_k)-1)+
\varphi_{m'_1,\ldots,m'_l}(i'_1,\ldots,i'_l)\\
&= \sum_{r=1}^k (i_r-1) n' n_{r+1} + \sum_{s=1}^l (i'_s-1) n'_{s+1} + 1\\
&=\varphi_{m_1,\ldots,m_k, m'_1,\ldots,m'_l} (i_1,\ldots,i_k,i'_1,\ldots,i'_l)
\end{align*}
where we have set $n_{k+1}=n'_{l+1}=1$.

\subsection[The union of the permutation groups]{} \label{UnionPermut}
We now consider the union
$$\hat{S} := \bigcup_{n=1}^\infty S_n$$
of all permutation groups\index{$\hat{S}$}\index{$S_n$}. Note that the different permutation groups are disjoint, which implies
that for an element $\sigma \in \hat{S}$ we can say exactly to which permutation group~$S_n$ it belongs.
This unique~$n$ will be called the degree of~$\sigma$.

On the set~$\hat{S}$, we introduce a product as follows: For $\sigma \in S_n$ and~$\tau \in S_m$, we define $\sigma \cdot \tau \in S_{mn}$\index{$\sigma \cdot \tau$} to be the unique permutation that makes the following diagram commutative:
$$
\xymatrix@C=6pc@R=4pc{
I_m \times I_n \ar[r]^{\varphi_{m,n}} \ar[d]_{\txt\scriptsize{$(i,j) \mapsto$ \\
$(\sigma(j),\tau(i))$}}& \ar[d]^{\txt\scriptsize{$\sigma \cdot \tau$}}
I_{mn} \\ I_n \times I_m \ar[r]_{\varphi_{n,m}}  & I_{mn} }
$$
Explicitly, this permutation is given by the formula
$$\sigma \cdot \tau ((i-1) n+j) = (\sigma(j)-1) m + \tau(i)$$
This product turns~$\hat{S}$ into a monoid:
\begin{prop}
$\hat{S}$ is a monoid with unit element~$\id_{I_1} \in S_1$.
\end{prop}
\begin{pf}
Suppose that $n$, $m$, and~$p$ are natural numbers. For $\rho \in S_p$, $\sigma \in S_n$, and $\tau \in S_m$, we have to show that
$(\rho \cdot \sigma) \cdot \tau = \rho \cdot (\sigma \cdot \tau)$. For this, note that the following diagram is commutative:
$$
\xymatrix@C=6pc@R=4pc{
I_m \times I_n \times I_p  \ar[r]^{\quad \varphi_{m,n} \times \id}
\ar[d]_{\txt\scriptsize{$(i,j,k) \mapsto$ \\ $(\rho(k),\sigma(j),\tau(i))$}}
& \ar[d]_{\txt\scriptsize{$(l,k) \mapsto$ \\ $(\rho(k),(\sigma \cdot\tau)(l))$}} I_{mn} \times
I_p \ar[r]^{\varphi_{mn,p}} & I_{mnp} \ar[d]^{\rho \cdot (\sigma \cdot \tau)} \\
I_p \times I_n \times I_m \ar[r]_{\quad \id \times \varphi_{n,m}}
& I_p \times I_{mn}  \ar[r]_{\varphi_{p,mn}} & I_{mnp}}
$$
In addition, the following diagram is also commutative:
$$
\xymatrix@C=6pc@R=4pc{
I_m \times I_n \times I_p  \ar[r]^{\quad  \id \times \varphi_{n,p}}
\ar[d]_{\txt\scriptsize{$(i,j,k) \mapsto$ \\ $(\rho(k),\sigma(j),\tau(i))$}}
& \ar[d]_{\txt\scriptsize{$(i,l) \mapsto$ \\ $((\rho \cdot \sigma)(l),\tau(i))$}}
I_m \times I_{np} \ar[r]^{\varphi_{m,np}} & I_{mnp} \ar[d]^{(\rho \cdot \sigma) \cdot \tau} \\
I_p \times I_n \times I_m \ar[r]_{\quad  \varphi_{p,n} \times \id}
& I_{np} \times I_m  \ar[r]_{\varphi_{np,m}} & I_{mnp}}
$$
Since the horizontal compositions in both diagrams are equal by Lemma~\ref{Monoton}, the associative law follows.

For the discussion of the unit element, note that the diagrams
$$
\xymatrix@C=6pc@R=4pc{
I_m \times I_1 \ar[r]^{\varphi_{m,1}} \ar[d]_{\txt\scriptsize{$(i,1) \mapsto$ \\
$(1,\tau(i))$}}& \ar[d]^{\tau}  I_m \\
I_1 \times I_m \ar[r]_{\varphi_{1,m}}  & I_m }
\qquad
\xymatrix@C=6pc@R=4pc{
I_1 \times I_n \ar[r]^{\varphi_{1,n}} \ar[d]_{\txt\scriptsize{$(1,j) \mapsto$ \\
$(\sigma(j),1)$}}& \ar[d]^{\sigma} I_n \\
I_n \times I_1 \ar[r]_{\varphi_{n,1}}  & I_n }
$$
commute, since $\varphi_{n,1}(i,1) = i = \varphi_{1,n}(1,i)$. This shows 
that~$\id_{I_1} \cdot \tau= \tau$ and
$\sigma \cdot \id_{I_1} = \sigma$. 
\qed
\end{pf}

\subsection[Bialgebras]{} \label{Bialg}
Suppose now that $H$ is a bialgebra with coproduct~$\Delta$ and counit~$\varepsilon$. We define certain powers that will play an important role in the sequel:
\begin{defn}
For $\sigma \in S_n$ and $h \in H$, we define the $\sigma$-th Sweedler power of~$h$ to be\index{$h^\sigma$}
$$h^\sigma := 
h_{(\sigma(1))} \cdot h_{(\sigma(2))} \cdot \ldots \cdot h_{(\sigma(n))}$$
\end{defn}
Here, we have used a variant of the so-called Sweedler notation 
(cf.~\cite{Sw}, Sec.~1.2, p.~10) to denote the images under the iterated comultiplication
$\Delta^n: H \rightarrow H^{\otimes n}$\index{$\Delta^n$} as 
$$\Delta^n(h) = h_{(1)} \otimes h_{(2)} \otimes \ldots \otimes h_{(n)}$$
If $n$ is the degree of~$\sigma$, we shall also say that $h^\sigma$ is an $n$-th Sweedler power of~$h$. 

This definition of Sweedler powers deviates slightly from the one in \cite{KashAnti}, \cite{KashPow}, where the permutation~$\sigma$ is always the identity. Of course, this permutation is not relevant if~$H$ is commutative or cocommutative, which is the case in the theory of algebraic groups, the setting in which Sweedler powers were first considered 
(cf.~\cite{SGA}, Par.~8.5, p.~474; \cite{TO}, Sec.~1, p.~3). But the introduction of the additional permutation causes Sweedler powers to be closed under iteration, as the
following power law for Sweedler powers asserts:
\begin{prop}
For $\sigma \in S_n$, $\tau \in S_m$, and~$h \in H$, we have
$(h^\sigma)^\tau = h^{\sigma \cdot \tau}$.
\end{prop}
\begin{pf}
We give a simplification of the original argument that was pointed out by P.~Schauenburg.
Recall that $\sigma \cdot \tau$ satisfies by definition the equation
$$\sigma \cdot \tau ((i-1) n+j) = (\sigma(j)-1) m + \tau(i)$$
We have in general that
$$\bigotimes_{j=1}^n \Delta^{m}(h_{(j)}) = h_\1 \o \ldots \o h_{(mn)} =
\bigotimes_{j=1}^n \bigotimes_{i=1}^m h_{((j-1)m+i)}$$
We therefore get
\begin{align*}
\Delta^{m}(h^\sigma) &= \Delta^{m}(\prod_{j=1}^n h_{(\sigma(j))}) = \prod_{j=1}^n \Delta^{m}(h_{(\sigma(j))})
= \prod_{j=1}^n \bigotimes_{i=1}^m h_{((\sigma(j)-1)m+i)} \\
&= \bigotimes_{i=1}^m \prod_{j=1}^n  h_{((\sigma(j)-1)m+i)}
\end{align*}
Permuting the tensor factors and multiplying, we get
\begin{align*}
(h^\sigma)^\tau = \prod_{i=1}^m \prod_{j=1}^n  h_{((\sigma(j)-1)m+\tau(i))}
= \prod_{i=1}^m \prod_{j=1}^n  h_{(\sigma \cdot \tau ((i-1) n+j))}
= \prod_{i=1}^{mn} h_{(\sigma \cdot \tau(i))} = h^{\sigma \cdot \tau}
\end{align*}
as asserted.
\qed
\end{pf}

\subsection[A monoid]{} \label{Monoid}
We proceed to construct certain Sweedler powers from special permutations that are based on finite sequences of positive integers. We denote such sequences, using square brackets, in the form $[n_1,n_2,\ldots,n_k]$. We call such a sequence normalized if every entry divides its predecessor, so that we have
$n_k / n_{k-1} / \ldots / n_2 / n_1$. Given a
sequence~$[n_1,n_2,\ldots,n_k]$, we define its normalization~$[n'_1,n'_2,\ldots,n'_k]$ recursively as follows: We set $n'_1 := n_1$ and 
$n'_{i+1} := \gcd(n_{i+1}, n'_i)$\index{$\gcd$}, the greatest common divisor of $n_{i+1}$ and $n'_i$.  

We define a product of two such sequences by the formula
$$[n_1,n_2,\ldots,n_k]  [m_1,m_2,\ldots,m_l] =
[m_1 n_1, m_1 n_2, \ldots, m_1 n_k, m_1, m_2, \ldots, m_l]$$
In addition, we introduce a unique element, called the empty sequence and denoted by~$[]$, that is by definition normalized and a unit for this product. If we denote the set of all such sequences by~$M$ and the set of all normalized sequences by~$M_N$, we have the following result:
\begin{prop}
$M$ is a monoid and~$M_N$ is a submonoid. Normalization 
defines a monoid homomorphism from $M$ to~$M_N$.
\end{prop}
\begin{pf}
The product is associative since we have
\begin{align*}
&([n_1,n_2,\ldots,n_k] [m_1,m_2,\ldots,m_l]) [p_1,p_2,\ldots,p_q] =\\
&[p_1 m_1 n_1, p_1 m_1 n_2, \ldots, p_1 m_1 n_k, p_1 m_1,
p_1 m_2, \ldots, p_1 m_l,p_1,p_2,\ldots,p_q] = \\
&[n_1,n_2,\ldots,n_k] ([m_1,m_2,\ldots,m_l] [p_1,p_2,\ldots,p_q])
\end{align*}
It is obvious that~$M_N$ is a submonoid. If $[n_1,\ldots,n_k]$ and $[m_1,\ldots,m_l]$ are two
sequences with normalizations $[n'_1,\ldots,n'_k]$ and $[m'_1,\ldots,m'_l]$, then the normalization
of the product sequence $[m_1 n_1, m_1 n_2, \ldots, m_1 n_k, m_1, m_2, \ldots, m_l]$ is
$$[m'_1 n'_1, m'_1 n'_2, \ldots, m'_1 n'_k, m'_1, m'_2, \ldots, m'_l]$$
since $m'_1 n'_{i+1} = \gcd(m'_1 n_{i+1}, m'_1 n'_i)$. This proves that normalization is a monoid homomorphism.
\qed
\end{pf}

\subsection[Permutations from sequences]{} \label{PermSeq}
Suppose now that $[n_1,n_2,\ldots,n_k]$ is a finite sequence of positive integers, and denote its normalization by $[n'_1,n'_2,\ldots,n'_k]$, so that $n'_{i+1} := \gcd(n_{i+1}, n'_i)$. We define the numbers $m_1,m_2,\ldots,m_{k-1}$ by $n'_i = m_i n'_{i+1}$ and put $m_k:=n'_k$. Similarly, we define $l_1,l_2,\ldots,l_{k-1}$ by~$n_{i+1}= l_i n'_{i+1}$ and put $l_k:=1$. Note that $m_i$ and~$l_i$ are relatively prime by construction, which implies that the map
$$\rho_i: I_{m_i} \rightarrow I_{m_i},~j \mapsto l_i (j-1) + 1$$
is bijective, where the right hand side really should be understood as the unique element of~$I_{m_i}$ that is congruent to~$l_i(j-1)+1$ modulo~$m_i$. 
The subtraction and addition of~$1$ in the definition of~$\rho_i$ stems from the fact that we have chosen the set~$I_{m_i}$ to consist of the numbers from $1$ to~$m_i$, and not from~$0$ to~$m_i-1$. Up to this shift, $\rho_i$ is just multiplication by~$l_i$. Note that $\rho_k=\id_{I_{m_k}}$.
\begin{defn}
We define $P(n_1,\ldots,n_k) \in S_{n_1}$\index{$P(n_1,\ldots,n_k)$} to be the unique
permutation that makes the following diagram commutative:
$$
\xymatrix@C=7pc@R=4pc{
I_{m_k} \times \ldots \times I_{m_1} \ar[r]^{\varphi_{m_k,\ldots,m_1}}
\ar[d]_{\txt\scriptsize{$(i_k,\ldots,i_1) \mapsto$ \\ $(\rho_1(i_1),\ldots,\rho_k(i_k))$}}
& I_{n_1} \ar[d]^{P(n_1,\ldots,n_k)} \\
I_{m_1} \times \ldots \times I_{m_k} \ar[r]_{\varphi_{m_1,\ldots,m_k}}
& I_{n_1}}
$$
For the empty sequence~$[]$, we define $P():=\id_{I_1}$.
\end{defn}
Since we have $n_1=m_1 \cdot \ldots \cdot m_k$, this definition makes sense. The relation between sequences and permutations becomes clear through the following fact:
\begin{prop}
The map
$$M \rightarrow \hat{S},~[n_1,\ldots,n_k] \mapsto P(n_1,\ldots,n_k)$$
is a monoid homomorphism.
\end{prop}
\begin{pf}
Suppose that, besides the sequence $[n_1,n_2,\ldots,n_k]$ considered above, we have another sequence $[p_1,p_2,\ldots,p_l]$ with normalization
$[p'_1,p'_2,\ldots,p'_l]$, so that $p'_{i+1}=\gcd(p_{i+1},p'_i)$. For $i=1,\ldots,l-1$, define $q_i$
by $p'_i = q_i p'_{i+1}$ and~$r_i$ by $p_{i+1} = r_i p'_{i+1}$. In addition, we set $q_l:=p'_l$ and $r_l:=1$. For $i=1,\ldots,l$, we then have a permutation
$$\pi_i : I_{q_i} \rightarrow I_{q_i},~j \mapsto r_i (j-1) + 1$$
and $P(p_1,\ldots,p_l)$ is defined via the commutative diagram
$$
\xymatrix@C=7pc@R=4pc{
I_{q_l} \times \ldots \times I_{q_1} \ar[r]^{\varphi_{q_l,\ldots,q_1}}
\ar[d]_{\txt\scriptsize{$(j_l,\ldots,j_1) \mapsto$ \\ $(\pi_1(j_1),\ldots,\pi_l(j_l))$}}
& I_{p_1} \ar[d]^{P(p_1,\ldots,p_l)} \\
I_{q_1} \times \ldots \times I_{q_l} \ar[r]_{\varphi_{q_1,\ldots,q_l}}
& I_{p_1}}
$$
Consider now the sequence
$$[n_1,n_2,\ldots,n_k][p_1,p_2,\ldots,p_l]=[n_1 p_1,n_2 p_1,\ldots,n_k p_1,p_1,p_2,\ldots,p_l]$$
By Proposition~\ref{Monoid}, its normalization is 
$[n'_1 p'_1,n'_2 p'_1,\ldots,n'_k p'_1,p'_1,p'_2,\ldots,p'_l]$. \linebreak[4]
Therefore, forming quotients of adjacent elements as above, we arrive at the sequence
$$[m_1,m_2,\ldots,m_k,q_1,q_2,\ldots,q_l]$$
By dividing every element of the sequence by its normalized counterpart and shifting the result by~1, we arrive at the sequence
$$[l_1,l_2,\ldots,l_{k-1},l_k,r_1,r_2,\ldots,r_{l-1},r_l]$$
We then have the commutative diagram
$$
\xymatrix@C=2.7pc@R=5pc{
I_{q_l} \times \ldots \times I_{q_1} \times
I_{m_k} \times \ldots \times I_{m_1} 
\ar[rr]^{\mspace{110mu} \varphi_{q_l,\ldots,q_1} \times \varphi_{m_k,\ldots,m_1}}
\ar[d]_{\txt\scriptsize{$(j_l,\ldots,j_1,i_k,\ldots,i_1) \mapsto$ \\
$(\rho_1(i_1),\ldots,\rho_k(i_k),$\\
$\mspace{40mu} \pi_1(j_1),\ldots,\pi_l(j_l))$}}
& &I_{p_1} \times I_{n_1} \ar[d]_{\txt\scriptsize{$(j,i) \mapsto$ \\ $(\sigma(i),\sigma'(j))$}}
\ar[r]^{\varphi_{p_1,n_1}} & I_{n_1 p_1}  \ar[d]^{\sigma''} \\
I_{m_1} \times \ldots \times I_{m_k} \times
I_{q_1} \times \ldots \times I_{q_l}
\ar[rr]_{\mspace{110mu} \varphi_{m_1,\ldots,m_k} \times \varphi_{q_1,\ldots,q_l}}
& & I_{n_1} \times I_{p_1} \ar[r]_{\varphi_{n_1,p_1}}
& I_{n_1 p_1} }
$$
where we have used the abbreviations
$$\sigma:=P(n_1,\ldots,n_k) \quad \sigma':=P(p_1,\ldots,p_l) \quad
\sigma'':=P(n_1p_1,\ldots,n_kp_1,p_1,\ldots,p_l)$$
In this diagram, the left rectangle is commutative by definition, and the large rectangle is commutative since the horizontal arrows are, by Lemma~\ref{Monoton}, equal to $\varphi_{q_l,\ldots,q_1,m_k,\ldots,m_1}$ and $\varphi_{m_1,\ldots,m_k,q_1,\ldots,q_l}$ respectively. Therefore, the right rectangle is also commutative, which, by comparison with the definition in Paragraph~\ref{UnionPermut}, shows that
$\sigma''= \sigma \cdot \sigma'$.
\qed
\end{pf}

\subsection[Sweedler powers]{} \label{SweedlPow}
Suppose now again that $H$ is a bialgebra and that $[n_1,n_2,\ldots,n_k]$ is a finite sequence of positive integers. For $h \in H$, we define\index{$h^{[n_1,n_2,\ldots,n_k]}$}
$$h^{[n_1,n_2,\ldots,n_k]} := h^{P(n_1,n_2,\ldots,n_k)}$$
and call it the
$[n_1,n_2,\ldots,n_k]$-th Sweedler power of~$h$. It is an immediate consequence of Proposition~\ref{Bialg} and Proposition~\ref{PermSeq} that we have the power law
$$(h^{[n_1,\ldots,n_k]})^{[p_1,\ldots,p_l]} = h^{[n_1,\ldots,n_k][p_1,\ldots,p_l]}$$
In the sequel, two special cases will be of particular importance for us. 
First, for sequences of length~1, we obviously have
$h^{[n]}= h_\1 \cdot h_\2 \cdot \ldots \cdot h_{(n)}$,
which establishes the connection with the notation used in~\cite{LinMont1}.
Second, consider sequences~$[n,k]$ of length~2 in which~$n$ and~$k$ are relatively prime. The corresponding normalized sequence then is~$[n,1]$. We get 
$$P(n,k)(1+j)=1+kj$$
and therefore
$$h^{[n,k]} = 
h_\1 \cdot h_{(1+k)} \cdot h_{(1+2k)} \cdot \ldots \cdot h_{(1+(n-1)k)}$$
where in the last two formulas the numbers on the right hand side have to be interpreted modulo~$n$, as explained in Paragraph~\ref{PermSeq}. For example, we have $h^{[5,2]} = h_{(1)}h_{(3)}h_{(5)}h_{(2)}h_{(4)}$. Of course, the second case reduces to the first if $k=1$, i.e., we have $h^{[n,1]}=h^{[n]}$.

\section{Frobenius-Schur Indicators}  \label{Sec:FrobSchur}
\subsection[Central Sweedler powers]{} \label{SweedlPowCent}
Recall that, as stated in the introduction, we assume from now on that~$K$ is an algebraically closed field of characteristic zero and that $H$ is a semisimple Hopf algebra over~$K$ with an integral~$\Lambda$ satisfying~$\varepsilon(\Lambda)=1$. We consider the Sweedler powers of this integral. The first interesting fact about them is that Sweedler powers of the kind just considered are central:
\begin{prop}
Suppose that $m$ and~$k$ are relatively prime natural numbers. Then~$\Lambda^{[m,k]}$ is central.
\end{prop}
\begin{pf}
From~\cite{LR2}, Lem.~1.2, p.~270, we have that 
$h\Lambda_\1 \o \Lambda_\2 = \Lambda_\1 \o S(h)\Lambda_\2$, which implies $$h\Lambda_\1 \o \Lambda_\2 \o \ldots \o \Lambda_{(m)}
= \Lambda_\1 \o \bigotimes_{i=2}^{m} S(h_{(m+1-i)}) \Lambda_{(i)}$$
To be able to read all Sweedler indices modulo~$m$, we write this in the form
$$h\Lambda_\1 \o \Lambda_\2 \o \ldots \o \Lambda_{(m)}
= \varepsilon(S(h_\1)) \Lambda_\1 \o \bigotimes_{i=2}^{m} S(h_{(m+2-i)}) \Lambda_{(i)}$$
By permuting the tensor factors and multiplying them afterwards, we arrive at
\begin{align*} h \Lambda^{[m,k]} &=
\varepsilon(S(h_\1)) \Lambda_\1 \prod_{i=1}^{m-1} S(h_{(m+2-(1+ik))}) \Lambda_{(1+ik)} \\
&= \varepsilon(S(h_\1)) \Lambda_\1 \prod_{i=1}^{m-1} S(h_{(1-ik)}) \Lambda_{(1+ik)}
\end{align*}
On the other hand, we have
$\Lambda_\1 S(h) \o \Lambda_\2 = \Lambda_\1 \o \Lambda_\2 h$, which implies
$$\Lambda_\1 \o \Lambda_\2 \o \ldots \o \Lambda_{(m)} h
= (\bigotimes_{i=1}^{m-1} \Lambda_{(i)} S(h_{(m-i)})) \o \Lambda_{(m)}$$
As above, we write this in the form
$$\Lambda_\1 \o \Lambda_\2 \o \ldots \o \Lambda_{(m)} h
= (\bigotimes_{i=1}^{m-1} \Lambda_{(i)} S(h_{(m+1-i)})) \o \Lambda_{(m)} \varepsilon( S(h_\1))$$
to be able to read the Sweedler indices modulo~$m$. Since $H$ is involutory, we know from~\cite{LR1}, Prop.~2(b), p.~191 that~$\Lambda$ is cocommutative, i.e., that 
$\Lambda_\1 \o \Lambda_\2 = \Lambda_\2 \o \Lambda_\1$, and therefore we can permute the Sweedler components of the integral in this expression cyclicly:
\begin{align*}
\bigotimes_{i=1}^{m-1} \Lambda&_{(i+1+(m-1)k)} \o \Lambda_{(1+(m-1)k)} h \\
&= (\bigotimes_{i=1}^{m-1} \Lambda_{(i+1+(m-1)k)} S(h_{(m+1-i)})) \o \Lambda_{(1+(m-1)k)} \varepsilon(S(h_\1))
\end{align*}
Note that the terms $\Lambda_{(i)} S(h_{(j)})$ appearing in the last formula all satisfy the requirement that $i+j \equiv 2-k \pmod{m}$. We now permute the first tensor factors and multiply afterwards. We then get
\begin{align*}
\Lambda^{[m,k]} h = (\prod_{i=0}^{m-1} \Lambda_{(1+ik)}) h
= (\prod_{i=0}^{m-2} \Lambda_{(1+ik)} S(h_{(1-(i+1)k)}))
\Lambda_{(1+(m-1)k)} \varepsilon(S(h_\1))
\end{align*}
By comparing both expressions, we get $h\Lambda^{[m,k]} =  \Lambda^{[m,k]}h$, as asserted.
\qed
\end{pf}
However, it is not true that arbitrary Sweedler powers of the integral are central, as we will see in Paragraph~\ref{SweedlPowNonCent} by an explicit example.

\subsection[The coproduct of the Sweedler powers]{} \label{CoprodSweedl}
As we will see even earlier in Paragraph~\ref{SweedlPowNonCocom}, the Sweedler powers of the integral are, in general, not cocommutative. However, there is a slightly weaker property that holds not only for the integral, but for an arbitrary cocommutative element:
\begin{prop}
Suppose that $m$ and~$k$ are relatively prime natural numbers, and that $x \in H$ and~$\chi \in H^*$ are cocommutative elements. Then we have
$$(\id \o \chi)(\Delta(x^{[m,k]})) = (\chi \o \id)(\Delta(x^{[m,k]}))$$
\end{prop}
\begin{pf}
Since $m$ and~$k$ are relatively prime, there exists a natural 
number~$l \le m$ such that $kl \equiv 1 \pmod{m}$. Let~$\sigma := P(m,k)$.
As in Paragraph~\ref{Bialg}, we have
$$\Delta(x^{[m,k]}) =
(\prod_{j=1}^m  x_{(2(\sigma(j)-1)+1)}) \o (\prod_{j=1}^m x_{(2(\sigma(j)-1)+2)})$$
As we have $\sigma(j)-1 \equiv k(j-1) \pmod{m}$, we have
$2(\sigma(j)-1) \equiv 2k(j-1) \pmod{2m}$, which gives
\begin{align*}
\Delta(x^{[m,k]})
&= (\prod_{j=1}^m  x_{(2k(j-1)+1)}) \o (\prod_{j=1}^m x_{(2k(j-1)+2)}) \\
&= (\prod_{j=0}^{m-1}  x_{(2kj+1)}) \o (\prod_{j=0}^{m-1} x_{(2kj+2)})
\end{align*}
where all Sweedler indices are read modulo~$2m$ to lie in~$I_{2m}$. We therefore get
\begin{align*}
(\chi& \o \id)(\Delta(x^{[m,k]}))
= \chi(\prod_{j=0}^{m-1}  x_{(2kj+1)}) \; \prod_{j=0}^{m-1} x_{(2kj+2)} \\
&= \chi(\prod_{j=0}^{m-1}  x_{(2kj)}) \; \prod_{j=0}^{m-1} x_{(2kj+1)} 
= \chi((\prod_{j=l}^{m-1}  x_{(2kj)})(\prod_{j=0}^{l-1}  x_{(2kj)})) \;
\prod_{j=0}^{m-1} x_{(2kj+1)} \\
&= \chi((\prod_{j=0}^{m-l-1}  x_{(2k(j+l))})(\prod_{j=m-l}^{m-1}  x_{(2k(j+l-m))})) \;
\prod_{j=0}^{m-1} x_{(2kj+1)} \\
&= \chi(\prod_{j=0}^{m-1}  x_{(2kj+2)}) \; \prod_{j=0}^{m-1} x_{(2kj+1)}
= (\id \o \chi)(\Delta(x^{[m,k]}))
\end{align*}
Here, the second equality follows from the cocommutativity of~$x$, the third equality follows from the cocommutativity of~$\chi$, and the fifth equality holds since $2kl \equiv 2 \pmod{2m}$.
\qed
\end{pf}

\subsection[Frobenius-Schur indicators]{} \label{FrobSchur}
Since the Sweedler powers $\Lambda^{[m]}$ are central by Proposition~\ref{SweedlPowCent}, they are determined by the values that the characters of~$H$ take on them. These values are the Frobenius-Schur indicators of the characters (cf.~\cite{LinMont1}, p.~348):
\begin{defn}
Suppose that~$m$ is a natural number and that~$\chi$ is a character of~$H$. The $m$-th Frobenius-Schur indicator of~$\chi$ is the number\index{$\nu_m(\chi)$}
$$\nu_m(\chi) := \chi(\Lambda^{[m]})$$
\end{defn}
The classical Frobenius-Schur indicator is the indicator~$\nu_2(\chi)$; it is understood that this indicator is meant when no~$m$ is explicitly mentioned. We will refer to the indicators~$\nu_m(\chi)$ for~$m>2$ as the higher indicators.

Let us consider the case of group rings. So, let~$G$\index{$G$} be a finite group and let $H=K[G]$\index{$K[G]$} be its group ring. The integral is then\index{$\vert G \vert$} 
$\Lambda={{\scriptstyle 1/}{}_{|G|}} \sum_{g \in G} g$. As the Sweedler power of a grouplike element coincides with its usual power, we have
$$\Lambda^{[m,k]} = \frac{1}{|G|} \sum_{g \in G} g^m$$
for the Sweedler powers of the integral. In particular, the right hand side of this expression is independent of~$k$. Therefore, the Frobenius-Schur indicators are in this case given by the formula
$$\nu_m(\chi) = \frac{1}{|G|} \sum_{g \in G} \chi(g^m)$$
(cf.~\cite{Isaacs}, Lem.~(4.4), p.~49).

Returning to the general case, we derive a first formula for the higher Frobenius-Schur indicators that should be viewed as a generalization of the Frobenius-Schur theorem for these indicators, since it reduces, as we will explain in Paragraph~\ref{FrobSchurThm}, in the case $m=2$ to the Frobenius-Schur theorem (cf.~\cite{LinMont1}, Thm.~3.1, p.~349). This formula follows from the following easy observation in linear algebra: Suppose that~$V$ is a finite-dimensional vector space and that $f_1,f_2,\ldots,f_m$ are linear endomorphisms of~$V$. Consider the linear map
$$\alpha: V^{\o m} \ra V^{\o m},~v_1 \o \ldots \o v_m \mapsto v_2 \o v_3 \o \ldots \o v_m \o v_1$$
About this simple situation, one can state the following general fact:
\begin{lemma}
Suppose that $k$ is relatively prime to~$m$, and let $\sigma:=P(m,k)$.
Then we have
$$\Tr_{V^{\o m}}(\alpha^k \circ (f_1 \o f_2 \o \ldots \o f_m))
=\Tr_V(f_{\sigma(1)} \circ f_{\sigma(2)} \circ \ldots \circ f_{\sigma(m)})$$
\end{lemma}
\begin{pf}
Suppose that $v_1,\ldots,v_n$ is a basis of~$V$ with dual basis~$v_1^*,\ldots,v_n^*$. For $l=1,\ldots,m$, we can represent~$f_l$ by a matrix that is determined by the formula
$$f_l(v_j) = \sum_{i=1}^n a_{ij}^l v_i$$
so that $a_{ij}^l= \langle v_i^*, f_l(v_j) \rangle$, where the angle brackets are used to denote the canonical pairing between a vector space and its dual. We now have
\begin{align*}
\Tr_{V^{\o m}}(&\alpha^k \circ (f_1 \o f_2 \o \ldots \o f_m)) \\
&= \sum_{i_1,\ldots,i_m=1}^n
\langle v_{i_1}^* \o v_{i_2}^* \o \ldots \o v_{i_m}^*,
\alpha^k(f_1(v_{i_1}) \o f_2(v_{i_2}) \o \ldots \o f_m(v_{i_m})) \rangle \\
&= \sum_{i_1,\ldots,i_m=1}^n
\langle v_{i_1}^*, f_{1+k}(v_{i_{1+k}}) \rangle
\langle v_{i_2}^*, f_{2+k}(v_{i_{2+k}}) \rangle  \ldots
\langle v_{i_m}^*, f_{m+k}(v_{i_{m+k}}) \rangle\\
&= \sum_{i_1,\ldots,i_m=1}^n  a_{i_1,i_{1+k}}^{1+k} \cdot a_{i_2,i_{2+k}}^{2+k} \cdot \ldots \cdot a_{i_m,i_{m+k}}^{m+k}
\end{align*}
where we read the indices modulo~$m$. We can now permute the factors in the product with the help of~$\sigma$, taking into account that
$\sigma(j+1)=\sigma(j)+k$:
$$\prod_{j=1}^m a_{i_j,i_{j+k}}^{j+k} = 
\prod_{j=1}^m a_{i_{\sigma(j)},i_{\sigma(j)+k}}^{\sigma(j)+k}= 
\prod_{j=1}^m a_{i_{\sigma(j)},i_{\sigma(j+1)}}^{\sigma(j+1)}$$
By replacing $i_1,\ldots,i_m$ by $i_{\sigma^{-1}(1)},\ldots,i_{\sigma^{-1}(m)}$ in the above summation, we get
\begin{align*}
\Tr_{V^{\o m}}(&\alpha^k \circ (f_1 \o f_2 \o \ldots \o f_m)) = \sum_{i_1,\ldots,i_m=1}^n \prod_{j=1}^m a_{i_j,i_{j+1}}^{\sigma(j+1)} \\
&= \Tr_V(f_{\sigma(2)} \circ f_{\sigma(3)} \circ \ldots \circ f_{\sigma(m)} \circ f_{\sigma(1)})
\end{align*}
The assertion now follows from the cyclicity of the trace.
\qed
\end{pf}
In particular, we get for $k=1$ that 
$$\Tr_{V^{\o m}}(\alpha \circ (f_1 \o f_2 \o \ldots \o f_m))
=\Tr_V(f_1 \circ f_2 \circ \ldots \circ f_m)$$
Now suppose that~$V$ is an $H$-module with character~$\chi$ and that 
\mbox{$\rho: H \rightarrow \End(V)$}\index{$\rho$} is the corresponding representation. From the diagonal action of~$H$ on the 
\mbox{$m$-th} tensor power, we then get the representation
$\rho^m: H \rightarrow \End(V^{\o m})$\index{$\rho^m$}.
The preceding lemma now relates the action of~$h \in H$ on the $m$-th tensor power directly to the Sweedler powers of~$h$:
\begin{prop}
Suppose that $k$ is relatively prime to~$m$.
Then we have
$$\Tr_{V^{\o m}}(\alpha^k \circ \rho^m(h))
= \chi(h^{[m,k]})$$
\end{prop}
If $h \in H$ is cocommutative, then the action of~$h$ on~$V^{\o m}$ commutes with~$\alpha$, since the equation
$h_\1 \o h_\2 = h_\2 \o h_\1$ implies that
$$h_\1 \o h_\2 \o \ldots \o h_\m= h_\2 \o h_\3 \o \ldots \o h_\m \o h_\1$$
This applies in particular to the integral~$\Lambda$, since the integral is, as discussed in Paragraph~\ref{SweedlPowCent}, a cocommutative element. As the action of~$\Lambda$ on~$V^{\o m}$ is the projection onto the subspace of invariants~$(V^{\o m})^H$, this space is invariant under~$\alpha$, and the trace of the restriction is the trace of the composition of~$\alpha$ and the projection. Therefore, the above proposition yields the following first formula for the Frobenius-Schur indicators:
\begin{corollary}
$\nu_m(\chi) = \Tr(\alpha \mid_{(V^{\o m})^H})$
\end{corollary}
It should be noted that, in general, the space $(V^{\o m})^H$ is not invariant under arbitrary permutations of the tensor factors. This corresponds to the fact that the trace of a product of matrices is only invariant under cyclic permutations of the matrices, but not under arbitrary permutations.

We will prove a second formula for the Frobenius-Schur indicators in Paragraph~\ref{SecondForm}. Note that the first formula implies that~$\nu_m(\chi)$ is a sum of certain eigenvalues of~$\alpha$, which are $m$-th roots of unity. It is therefore contained 
in~$\Q_m \subset K$\index{$\Q_m$}, the $m$-th cyclotomic field.

It is clear that the dual module~$V^*$ has the same Frobenius-Schur indicators as~$V$: Since~$S(\Lambda) = \Lambda$, we have
$$S(\Lambda^{[m]})= S(\Lambda_\1 \Lambda_\2 \ldots \Lambda_\m)) 
= S(\Lambda_\m) \ldots S(\Lambda_\2) S(\Lambda_\1) = \Lambda^{[m]}$$
and therefore we have for the character~$\bar{\chi}=S(\chi)$\index{$\bar{\chi}$} of~$V^*$ that $\bar{\chi}(\Lambda^{[m]}) = \chi(\Lambda^{[m]})$. However, as we will see in Paragraph~\ref{ExampleA4}, this does not imply that the higher indicators are always real numbers.

\subsection[The Frobenius-Schur theorem]{} \label{FrobSchurThm}
In \cite{FS}, F.~G.~Frobenius and I.~Schur investigated the question whether a given simple module of a finite group~$G$ over the complex numbers admits a basis such that the corresponding matrix representations of the actions of the group elements have real entries. They realized that the  necessary condition that the character~$\chi$ takes real values on all group elements is not sufficient. To find a sufficient condition, they introduced what is
now called the Frobenius-Schur indicator, and proved that a sufficient condition is that the indicator of this character is~$1$.

Since $\overline{\chi(g)}=\chi(g^{-1})$, it is clear that the condition that the character takes only real values is equivalent to the self-duality of the representation. It is not so clear that the matrix representations can be realized by real matrices if and only if there is a nondegenerate symmetric invariant bilinear form on the module (cf.~\cite{Serre1}, Sec.~13.2, Thm.~31, p.~106 for a proof). Using this, however, the result can be carried over to other base fields, and even to Hopf algebras. This was realized by V.~Linchenko and S.~Montgomery (cf.~\cite{LinMont1}), and also in similar form by J.~Fuchs, A.~Ch. Ganchev, K.~Szlach\'{a}nyi, and P.~Vecserny\'{e}s (cf.~\cite{FuGanSzlVec1}). Their Frobenius-Schur theorem for Hopf algebras states, under the assumptions made in this section, the following (cf.~\cite{LinMont1}, Thm.~3.1, p.~349; see also \cite{FuGanSzlVec1}, Sec.~IV):
\begin{thm}
Suppose that~$V$ is a simple $H$-module with character~$\chi$. Then its Frobenius-Schur indicator is
$$\nu_2(\chi) = 
\begin{cases}
0:& V \; \text{is not self-dual} \\
1:& V \; \text{admits a nondegenerate} \\
& \text{symmetric invariant bilinear form} \\
-1:& V \; \text{admits a nondegenerate skew-} \\
&  \text{symmetric invariant bilinear form} 
\end{cases}$$
\end{thm}
\begin{pf}
We derive this from Corollary~\ref{FrobSchur}. Recall first that a bilinear form 
$\langle \cdot,\cdot \rangle: V \times V \rightarrow K$ 
is called invariant if
$$\langle h_\1.v, h_\2.v' \rangle = \varepsilon(h) \langle v, v' \rangle$$
for all $v, v' \in V$ and all~$h \in H$, which is of course equivalent to the condition that 
$\langle h.v, v' \rangle = \langle v, S(h).v' \rangle$. We have the isomorphism
$$(V \o V)^H \rightarrow \Hom_H(V^*,V),~x \mapsto 
(\varphi \mapsto (\id \o \, \varphi) x)$$
Schur's lemma therefore implies that $(V \o V)^H \neq 0$ if and only 
if~$V \cong V^*$. In particular, if~$V$ is not self-dual, then Corollary~\ref{FrobSchur} yields
$$\nu_2(\chi) = \Tr(\alpha \mid_{(V \o V)^H}) = 0$$
establishing the first case of the Frobenius-Schur theorem.

Suppose now that~$V$ is self-dual. As we have just seen, homomorphisms
from~$V$ to~$V^*$ are given by invariant bilinear forms. Since the space of homomorphisms is one-dimensional, an invariant bilinear form is unique up to scalar multiples, and since every nonzero homomorphism is an isomorphism, a nonzero invariant bilinear form is nondegenerate. If 
$\langle \cdot,\cdot \rangle$ is a nondegenerate invariant bilinear form, the bilinear form that arises by interchanging the arguments is also invariant, as we have
$$\langle S(h).v', v \rangle = \langle v', S^2(h).v \rangle
= \langle v', h.v \rangle$$
Of course, we have crucially used here the fact that~$H$ is involutory. This implies that the form with interchanged arguments is proportional to the original form, so that
$$\langle v', v \rangle = \mu \langle v, v' \rangle$$
for a scalar~$\mu \in K$\index{$\mu$}. Interchanging the arguments another time, we arrive at the equation
$\langle v, v' \rangle = \mu^2 \langle v, v' \rangle$,
which yields $\mu = \pm 1$. We therefore see that our bilinear form is either symmetric or skew-symmetric.

Our invariant bilinear form determines a
tensor~$y = \sum_i \varphi_i \o \varphi'_i \in V^* \o V^*$\index{$y$} by the requirement that
$$\sum_i \varphi_i(v) \varphi'_i(v') = \langle v, v' \rangle$$
for all $v,v' \in V$. The invariance of the bilinear form translates to the condition
$\sum_i h_\1.\varphi_i \o h_\2.\varphi'_i 
= \varepsilon(h) \sum_i \varphi_i \o \varphi'_i$, which says that~$y$ is a nonzero element in the one-dimensional space of invariants 
$(V^* \o V^*)^H$. Furthermore\index{$V^H$}, the equation
$\langle v', v \rangle = \mu \langle v, v' \rangle$ yields 
$\alpha(y) = \mu y$. But now Corollary~\ref{FrobSchur} gives
$$\nu_2(\chi) = \Tr(\alpha \mid_{(V \o V)^H}) 
= \Tr(\alpha \mid_{(V^* \o V^*)^H})= \mu$$
which establishes the remaining two cases of the Frobenius-Schur theorem.
\qed
\end{pf}

\subsection[The Frobenius-Schur indicators of the regular representation]{} \label{FrobSchurRegRep}
It is in general difficult to compute the Frobenius-Schur indicators of a given module. However, in the case of the regular representation, in which $H$ acts on itself via left multiplication, there is a closed formula:
\begin{prop}
Suppose that $\chi_R$\index{$\chi_R$} is the character of the regular representation and consider the map
$$E_m: H \rightarrow H,~h \mapsto S(h^{[m-1]})$$
for a natural number~$m \ge 2$. Then the $m$-th Frobenius-Schur indicator of~$\chi_R$ is $\nu_m(\chi_R) = \Tr(E_m)$.\index{$E_m$}
\end{prop}
\begin{pf}
If $V$ is any $H$-module, then the map
$$V \o H \rightarrow V \o H,~v \o h \mapsto h_\1.v \o h_\2$$
is an $H$-linear isomorphism if the left hand side is endowed with the $H$-module structure given by left multiplication on the $H$-tensorand and the right hand side carries the usual diagonal module structure. From this, it is clear that the map
$$f: V \rightarrow (V \o H)^H,~v \mapsto \Lambda_\1.v \o \Lambda_\2$$
is bijective. We apply this with $V=H^{\o (m-1)}$, so that
$$
f(h_1 \o h_2 \o \ldots \o h_{m-1}) =
\Lambda_\1 h_1 \o \Lambda_\2 h_2 \o \ldots \o \Lambda_{(m-1)} h_{m-1} \o \Lambda_\m$$
If we introduce the map
\begin{align*}
&\mspace{200mu}
\beta: H^{\o (m-1)} \rightarrow H^{\o (m-1)}, \\
&h_1 \o h_2 \o \ldots \o h_{m-1} \mapsto \\
&\mspace{100mu}
S(h_1{}_{(m-1)}) h_2 \o S(h_1{}_{(m-2)}) h_3 \o \ldots \o S(h_1{}_\2) h_{m-1} \o S(h_1{}_\1)
\end{align*}
then the diagram
$$
\xymatrix@C=7pc@R=4pc{
H^{\o (m-1)} \ar[r]^{\beta} \ar[d]_{f} & H^{\o (m-1)} \ar[d]^{f} \\
(H^{\o m})^H \ar[r]_{\alpha} & (H^{\o m})^H}
$$
is commutative. This holds since from the equation
$\Lambda_\1 \o \Lambda_\2 h = \Lambda_\1 S(h) \o \Lambda_\2$ (cf.~\cite{LR2}, Lem.~1.2, p.~270) we have
\begin{align*}
\Lambda_\1 & S(h_{(m-1)})  \o \Lambda_\2 S(h_{(m-2)})  \o \ldots
\o \Lambda_{(m-1)}  S(h_\1) \o \Lambda_\m \\
&\mspace{250mu}
= \Lambda_\1  \o \Lambda_\2  \o \ldots \o \Lambda_{(m-1)}  \o \Lambda_\m h
\end{align*}
and therefore 
\begin{align*}
\alpha(&f(h_1 \o \ldots \o h_{m-1}))  \\
&=\Lambda_\2 h_2 \o \Lambda_\3 h_3 \o \ldots \o \Lambda_{(m-1)} h_{m-1} \o \Lambda_\m \o \Lambda_\1 h_1 \\
&=\Lambda_\1 h_2 \o \Lambda_\2 h_3 \o \ldots \o \Lambda_{(m-2)} h_{m-1} \o \Lambda_{(m-1)} \o \Lambda_\m h_1 \\
&=\Lambda_\1 S(h_{1(m-1)}) h_2 \o \Lambda_\2 S(h_{1(m-2)}) h_3 \o \ldots \\
&\mspace{200mu}
\o \Lambda_{(m-2)} S(h_1{}_\2) h_{m-1} \o \Lambda_{(m-1)}S(h_1{}_\1) \o \Lambda_\m \\
&= f(\beta(h_1 \o \ldots \o h_{m-1}))
\end{align*}
where we have used for the second equality that $\Lambda$ is cocommutative (cf.~\cite{LR1}, Prop.~2(b), p.~191). This implies that $\Tr(\alpha \mid_{(H^{\o m})^H}) = \Tr(\beta)$. To calculate the trace of~$\beta$,
choose a basis $b_1,\ldots,b_n \in H$ with dual basis $b^*_1,\ldots,b^*_n \in H^*$. We then have
\allowdisplaybreaks{
\begin{align*}
&\Tr(\beta) =
\sum_{i_1,\ldots,i_{m-1}=1}^n \langle b^*_{i_1} \o \ldots \o b^*_{i_{m-1}},
\beta(b_{i_1} \o \ldots \o b_{i_{m-1}}) \rangle \\
&= \sum_{i_1,\ldots,i_{m-1}=1}^n
b^*_{i_1}(S(b_{i_1}{}_{(m-1)}) b_{i_2})  b^*_{i_2}(S(b_{i_1}{}_{(m-2)}) b_{i_3})  \ldots \\
&\mspace{250mu}
b^*_{i_{m-2}}(S(b_{i_1}{}_\2) b_{i_{m-1}}) b^*_{i_{m-1}}(S(b_{i_1}{}_\1))  \\
&= \sum_{i_1,\ldots,i_{m-2}=1}^n
b^*_{i_1}(S(b_{i_1}{}_{(m-1)}) b_{i_2})  b^*_{i_2}(S(b_{i_1}{}_{(m-2)}) b_{i_3})  \ldots
b^*_{i_{m-2}}(S(b_{i_1}{}_\2) S(b_{i_1}{}_\1))   \\
&= \sum_{i_1,\ldots,i_{m-3}=1}^n
b^*_{i_1}(S(b_{i_1}{}_{(m-1)}) b_{i_2})  b^*_{i_2}(S(b_{i_1}{}_{(m-2)}) b_{i_3})  \ldots\\
&\mspace{250mu}
b^*_{i_{m-3}}(S(b_{i_1}{}_\3) S(b_{i_1}{}_\2) S(b_{i_1}{}_\1))   \\
&= \ldots =
\sum_{i_1=1}^n
b^*_{i_1}(S(b_{i_1}{}_{(m-1)}) S(b_{i_1}{}_{(m-2)})   \ldots
S(b_{i_1}{}_\3) S(b_{i_1}{}_\2) S(b_{i_1}{}_\1))   \\
&= \;
\sum_{i=1}^n
b^*_i(S(b_i{}_\1 b_i{}_\2 b_i{}_\3 \ldots b_i{}_{(m-2)} b_i{}_{(m-1)}))
= \Tr(E_m)
\end{align*}}
as asserted.
\qed
\end{pf}
In the case where $H=K[G]$ is the group ring of a finite group~$G$, the map~$E_m$ is given on the elements~$g \in G$ as
$E_m(g) = g^{-m+1}$. This shows that~$E_m$ permutes this basis of~$H$, and therefore its trace is the number of fixed points. It is clear that~$g$ is a fixed point of~$E_m$ if and only if~$g^m=1$, and therefore we get
$$\nu_m(\chi_R) = \Tr(E_m) = |\{g \in G \mid g^m=1\}|$$
so that the $m$-th Frobenius-Schur indicator of the regular representation is just the number of elements whose order divides~$m$. Of course, this can also be seen from the general form of the character of the regular representation (cf.~\cite{Serre1}, Sec.~2.4, Prop.~5, p.~18).

\section{The Exponent}  \label{Sec:Exponent}
\subsection[The exponent]{} \label{Exponent}
A second approach to Frobenius-Schur indicators uses the so-called exponent of a module, which we define now:
\begin{defn}
Suppose that~$V$ is an $H$-module. The exponent of~$V$, denoted by~$\exp(V)$\index{$\exp(V)$}, is the smallest natural number~$m$ such that
$$h^{[m]}.v = \varepsilon(h) v $$
for all $h \in H$ and all~$v \in V$. If no such~$m$ exists, we say that
the exponent of~$V$ is infinite. The exponent of~$H$ is the exponent of the regular representation.
\end{defn}
Note that this definition has to be modified if~$H$ is not involutory 
(cf.~\cite{EG1}, Def.~2.1, p.~132). In the case where $H=K[G]$ is the group ring of a finite group~$G$, the above condition becomes 
$g^m.v = v$ for all $g \in G$ and all~$v \in V$, and therefore the exponent of~$V$ is the exponent of the image of~$G$ in the general linear group~$GL(V)$\index{$GL(V)$}.

To understand this notion better, we approach it from another point of view. Choose a basis~$b_1,\ldots,b_n$ of~$H$ with dual basis~$b^*_1,\ldots,b^*_n$. The tensor\index{$t$} 
$$t:= \sum_{i=1}^n b_i \o b^*_i$$
does not depend on the choice of the basis. It acts on~$V \o H^*$ in such a way that the $H$-tensorand acts by the module action and the $H^*$-tensorand acts by left multiplication. This obviously implies that this action is $H^*$-linear with respect to the right $H^*$-module structure on~$V \o H^*$ that arises from right multiplication on the second tensorand. We denote the corresponding endomorphism of~$V \o H^*$ by~$t\mid_{V \o H^*}$. It is related to the exponent in the following way:
\begin{prop} 
The order of~$t\mid_{V \o H^*}$ is the exponent of~$V$.
\end{prop}
\begin{pf}
The $m$-th power of~$t$ is
$$t^m = \sum_{i_1,i_2,\ldots,i_m=1}^n 
b_{i_1} b_{i_2} \ldots b_{i_m} \o b^*_{i_1} b^*_{i_2} \ldots b^*_{i_m}$$
Evaluating the $H^*$-tensorand on~$h \in H$, we get
$$\sum_{i_1,i_2,\ldots,i_m=1}^n b_{i_1} b_{i_2} \ldots b_{i_m}  b^*_{i_1}(h_\1) b^*_{i_2}(h_\2) \ldots b^*_{i_m}(h_\m) = h^{[m]}$$
The condition $t^m\mid_{V \o H^*} =(1 \o \varepsilon)\mid_{V \o H^*}$ is therefore equivalent to the condition that every Sweedler power~$h^{[m]}$
acts on~$V$ by multiplication with~$\varepsilon(h)$. This gives the assertion.
\qed
\end{pf}
The action of~$t$ can be interpreted in another way: The left $H$-module~$V$ gives rise to a right $H^*$-coaction by the formula
$$\delta(v) := \sum_{i=1}^n b_i.v \o b^*_i$$
(cf.~\cite{M}, Lem.~1.6.4, p.~11). The action of~$t$ on~$V \o H^*$ can therefore be written in the form
$$t.(v \o \varphi) = v^{\1} \o v^{\2} \varphi$$
where we have used the Sweedler notation $\delta(v)= v^{\1} \o v^{\2}$ for the coaction. If $v_1,\ldots,v_k$ is a basis of~$V$, then we define the 
$k \times k$-matrix~$(c_{ij})_{i,j=1,\ldots,k}$ with entries in~$H^*$ by the equation
$$\delta(v_j) = \sum_{i=1}^k v_i \o c_{ij}$$
The comodule condition then implies that we have
$$\Delta(c_{ij}) = \sum_{q=1}^k c_{iq} \o c_{qj} \qquad \qquad
\varepsilon(c_{ij}) = \delta_{ij}$$
which also shows that the elements~$c_{ij}$ span a subcoalgebra of~$H^*$.
Now the elements $v_1 \o \varepsilon,\ldots,v_k \o \varepsilon$ form a basis for the right $H^*$-module~$V \o H^*$ considered above, and the matrix representation of~$t\mid_{V \o H^*}$ with respect to this basis is exactly the matrix~$(c_{ij})_{i,j=1,\ldots,k}$:
$$t.(v_j \o \varepsilon) = v_j^{\1} \o v_j^{\2} = \sum_{i=1}^k (v_i \o \varepsilon)c_{ij}$$
In particular, we have the following:
\begin{corollary} 
The order of the $k \times k$-matrix~$(c_{ij})$ is equal to the exponent of~$V$.
\end{corollary}
In the case where~$V$ is one-dimensional, every $h \in H$ acts by multiplication with a scalar~$\gamma(h)$, and therefore the coaction is given by the formula $\delta(v)= v \o \gamma$. The above matrix then is a $1 \times 1$-matrix with unique entry~$\gamma$, and therefore the above corollary shows that the exponent of a one-dimensional module is just the
order of the corresponding character in the group~$G(H^*)$\index{$G(H^*)$} of grouplike elements of~$H^*$.

The notion of the exponent of a Hopf algebra was introduced by P.~Etingof and S.~Gelaki (cf.~\cite{EG1}, Def.~2.1, p.~132). However, the concept was, without this name, in use earlier. It was first considered in the theory of algebraic groups, i.e., for commutative Hopf algebras (cf.~\cite{SGA}, Prop.~8.5, p.~473; \cite{TO}, Sec.~1, p.~4), where is arises naturally by taking powers of group elements. In the case where the Hopf algebra is not commutative and not cocommutative, the concept has not such an easy interpretation, but it also arises naturally in the investigation of the order of the antipode of a Yetter-Drinfel'd Hopf algebra (cf.~\cite{KashAnti}, Thm.~6, p.~1264; see also~\cite{So5}, Par.~3.8, p.~443). This lead the first author to the conjecture that the exponent of a semisimple Hopf algebra divides its dimension (cf.~\cite{KashAnti}, p.~1261; see also \cite{KashPow}, p.~159). Although this conjecture is still open at present, P.~Etingof and S.~Gelaki have proved the following important result (cf.~\cite{EG1}, Thm.~4.3, p.~136): 
\begin{thm}
$\exp(H)$ divides $\dim(H)^3$.
\end{thm}
This result will be used extensively in the sequel.

\subsection[The second formula for the Frobenius-Schur indicators]{} \label{SecondForm}
The proof of Proposition~\ref{Exponent} leads directly to the second formula for the Frobenius-Schur indicators, which expresses the indicators in terms of the tensor~$t$:
\begin{prop}
Suppose that~$V$ is an $H$-module with character~$\chi$ and that~$m$ is a natural number. Then the $m$-th Frobenius-Schur indicator is
$$\nu_m(\chi) = \frac{1}{\dim(H)} \Tr(t^m\mid_{V \o H^*})$$
\end{prop}
\begin{pf}
Recall that the trace of the left multiplication by~$\varphi \in H^*$ on~$H^*$ is~$\dim(H) \varphi(\Lambda)$ (cf.~\cite{LR2}, Prop.~2.4, p.~273). Since the trace of the tensor product of two matrices is the product of the traces, we get
\begin{align*}
\Tr(t^m\mid_{V \o H^*}) 
&= \dim(H) \sum_{i_1,i_2,\ldots,i_m=1}^n 
\chi(b_{i_1} b_{i_2} \ldots b_{i_m}) 
(b^*_{i_1} b^*_{i_2} \ldots b^*_{i_m})(\Lambda) \\
&= \dim(H) \chi(\Lambda^{[m]})
\end{align*}
as asserted.
\qed
\end{pf}
We will give a third formula for the Frobenius-Schur indicators in Paragraph~\ref{ThirdForm}. This third formula should be considered as a variant of the one above, although its statement requires the Drinfel'd double construction, which we will discuss in Section~\ref{Sec:DrinfDouble}.

If $d:=\exp(V)$\index{$d$}, it is clear from the second formula that~$\nu_m(\chi)$ is contained in the cyclotomic field~$\Q_d \subset K$, since the eigenvalues of~$t$ 
on~$V \o H^*$ are $d$-th roots of unity, and the trace is the sum of these eigenvalues. Therefore, combining the first and the second formula for the Frobenius-Schur indicators, we see that the Frobenius-Schur indicator~$\nu_m(\chi)$ is an integer if~$m$ is large enough in the sense that it includes enough prime factors of~$d$. More precisely, if $p$ is a prime that divides~$d$ as well as~$m$, then we require that~$p$ appears in~$m$ at least as many times as it appears in~$d$, so that, if 
$$d = p^k q, \quad p \nmid q$$
then $p^k$ divides~$m$. If this condition holds for all primes~$p$, we say that~$m$ is large compared to~$d$. Note that this is always the case if~$d$ is squarefree.
\begin{corollary}
Suppose that $V$ is an $H$-module with character~$\chi$ and exponent~$d$, and that~$m$ is a natural number.
\begin{enumerate}
\item If~$m$ is large compared to~$d$, then~$\nu_m(\chi)$ is an integer.
\item
If~$d$ is squarefree, $\nu_m(\chi)$ is always an integer.
\end{enumerate} 
\end{corollary}
\begin{pf}
Let~$d'$ be the greatest common divisor of~$m$ and~$d$. The condition that~$m$ is large compared to~$d$ can be reformulated by saying that $d/d'$ is relatively prime to~$m$.
Namely, if $p$ were a prime that divides both $d/d'$ and~$m$, then it clearly would divide both~$d$ and~$m$. If $d = p^k q$, where $p \nmid q$, then $p^k$ would divide~$m$ by assumption, and therefore it would also divide~$d'$. This would imply that~$p$ does not divide~$d/d'$, a contradiction. The converse is also clear: If $p$ is a prime that divides~$d$ as well as~$m$, then write
$$d = p^k q, \; p \nmid q \qquad m = p^l q', \; p \nmid q'$$
If now $k$ were bigger than~$l$, the prime power $p^{k-l}$ would still divide~$d/d'$, contrary to our assumption.

Now write $m=d'm'$. Since $t^{d'}\mid_{V \o H^*}$ has order~$d/d'$, we get from the second formula for the Frobenius-Schur indicators that
$$\nu_m(\chi) = \frac{1}{\dim(H)} \Tr((t^{d'})^{m'} \mid_{V \o H^*})\in \Q_{d/d'}$$
On the other hand, we get from Corollary~\ref{FrobSchur}, the first formula for the Frobenius-Schur indicators, that $\nu_m(\chi) \in \Q_m$. Since~$d/d'$ and~$m$ are relatively prime, we see that 
$\nu_m(\chi) \in \Q_{d/d'} \cap \Q_m = \Q$, so that it is a rational number (cf.~\cite{Wash}, Prop.~2.4, p.~11). But by the first formula, it is also an algebraic integer, and therefore it must be an integer (cf.~\cite{SchejaStorch2}, \S~56, Bsp.~3, p.~91). This proves the first statement. The second statement is an immediate consequence of the first. 
\qed
\end{pf}
This result suggest the conjecture that the Frobenius-Schur indicators are in fact always integers. However, we will see in Paragraph~\ref{ExampleA4} that this is not the case. But it is the case for group rings, as we will see in the next paragraph.

\subsection[Sweedler powers of the integral]{} \label{SweedlPowInt}
The fact that the Frobenius-Schur indicators are contained in certain cyclotomic fields leads to some curious relations between Sweedler powers of the integral of different order. Choose a primitive $m$-th root of unity~$\zeta$ in~$K$. The subfield generated by the prime field~$\Q$ and~$\zeta$ is the $m$-th cyclotomic field~$\Q_m \subset K$; its isomorphism type does not depend on the field~$K$ in which the construction is carried out. Therefore, there is the standard isomorphism
$$\Z_m^\times \rightarrow \Gal(\Q_m/\Q),~\bar{k} \mapsto \sigma_k$$
between the group of units~$\Z_m^\times$ of~$\Z_m$\index{$\Z_n$}\index{$\Z_n^\times$} and the Galois group~$\Gal(\Q_m/\Q)$ of the cyclotomic field, where~$\sigma_k$\index{$\sigma_k$} is determined by the equation
$\sigma_k(\zeta) = \zeta^k$.
From the two formulas for the Frobenius-Schur indicators, we can determine how the Galois group acts on the indicators:
\begin{prop}
Suppose that~$V$ is an $H$-module with character~$\chi$ and exponent~$d$, and that $m$, $k$, and~$l$ are natural numbers.
\begin{enumerate}
\item If $k$ and~$l$ are relatively prime to~$m$, we have $\sigma_l(\chi(\Lambda^{[m,k]})) = \chi(\Lambda^{[m,kl]})$.
\item If $l$ is relatively prime to~$d$, we have $\sigma_l(\chi(\Lambda^{[m]})) = \chi(\Lambda^{[ml]})$.
\end{enumerate} 
\end{prop}
\begin{pf}
Recall from the discussion after Proposition~\ref{FrobSchur} the formula
$$\Tr(\alpha^k \mid_{(V^{\o m})^H})
= \chi(\Lambda^{[m,k]})$$
where $\alpha \in \End(V^{\o m})$ is defined by
$$\alpha(v_1 \o \ldots \o v_m) = v_2 \o v_3 \o \ldots \o v_m \o v_1$$
Since $\alpha$ has order~$m$, its eigenvalues, and also the eigenvalues of~$\alpha^k$, are $m$-th roots of unity. The automorphism $\sigma_l \in \Gal(\Q_m/\Q)$ raises these eigenvalues to the $l$-th power, so that the eigenvalues of~$\alpha^{k}$ become the eigenvalues of~$\alpha^{kl}$.
Since the trace is the sum of the eigenvalues, the first assertion follows.
The second assertion follows similarly from Proposition~\ref{SecondForm}, the second formula for the Frobenius-Schur indicators:
$$\sigma_l(\chi(\Lambda^{[m]})) = 
{\textstyle \frac{1}{\dim(H)}} \sigma_l(\Tr(t^m \mid_{V \o H^*})) = 
{\textstyle \frac{1}{\dim(H)}} \Tr(t^{ml} \mid_{V \o H^*}) = 
\chi(\Lambda^{[ml]})$$
Note that in this equation $\sigma_l$ is an element of the Galois group of the cyclotomic field~$\Q_d$, and not of~$\Q_m$ as before. However, both definitions are compatible if both are defined, since both raise roots of unity to their $l$-th power.
\qed
\end{pf}
The first formula above tells in the case $k=1$ that
$\sigma_l(\chi(\Lambda^{[m]})) = \chi(\Lambda^{[m,l]})$.
In the case where $H=K[G]$ is the group ring of a finite group~$G$, we have
$\Lambda^{[m]} = \Lambda^{[m,l]}$, since an iterated coproduct of an element of a group ring does not change under any permutation of the tensor factors. This implies that the Frobenius-Schur indicators are invariant under the action of the Galois group, which implies that they are rational. As they are always algebraic integers by Corollary~\ref{FrobSchur}, the first formula for the indicators, this proves the known fact that the higher Frobenius-Schur indicators are integers in the case of group rings (cf.~\cite{Isaacs}, Exerc.~(4.7), p.~60).

From the above formulas for the action of the Galois group, we can deduce the promised relations between the Sweedler powers of the integral: 
\begin{corollary}
Suppose that $m$ and~$l$ are relatively prime natural numbers.
\begin{enumerate}
\item If $l$ is also relatively prime to~$\exp(H)$, we have 
$\Lambda^{[m,l]}= \Lambda^{[ml]}$.
\item
If $m$ is large compared to~$\exp(H)$, we have
\mbox{$\Lambda^{[m]} = \Lambda^{[m,l]}$.}
\end{enumerate} 
\end{corollary}
\begin{pf}
Note first that the order of~$t$ is the order of the action of this element on~$H \o H^*$, and therefore equal to~$\exp(H)$. This shows that, for every $H$-module~$V$ with character~$\chi$, $d:=\exp(V)$ divides~$\exp(H)$. 

For the first assertion, consider the element~$\sigma_l$ in the Galois group of the cyclotomic field~$\Q_{md}$. Note that its restrictions to the cyclotomic fields~$\Q_m$ and~$\Q_d$ are also denoted by~$\sigma_l$. Now the preceding proposition tells that
$$\chi(\Lambda^{[m,l]})= \sigma_l(\chi(\Lambda^{[m]})) = \chi(\Lambda^{[ml]})$$
Since these Sweedler powers of the integral are central by Proposition~\ref{SweedlPowCent}, they are determined by the values that the characters take on them. This establishes the first assertion. The second assertion follows by a similar reasoning from Corollary~\ref{SecondForm}.
\qed
\end{pf}
Besides being true, this corollary is in fact quite surprising, as it asserts,  for example, that for an odd-dimensional semisimple Hopf algebra we have
$$\Lambda_\1 \Lambda_\2 \Lambda_\3 \Lambda_\4 \Lambda_\5 \Lambda_\6 =
\Lambda_\1 \Lambda_\3 \Lambda_\2$$
This is the case $m=3$, $l=2$ of the first assertion in the above corollary. 

In this context, it should be noted that the Sweedler powers of the integral
determine the exponent completely: The exponent of~$V$ is the smallest natural number~$m$ such that $\Lambda^{[m]}$ acts as the identity on~$V$. This holds since by Proposition~\ref{SecondForm}, the second formula for the Frobenius-Schur indicator, we then have that the trace of~$t^m$ on~$V \o H^*$ is $\dim(V) \dim(H)$. As the trace is the sum of the eigenvalues, this is a sum of roots of unity that is equal to the number of summands. This can only happen if all these roots of unity are equal to one, so that we get that 
$t^m \mid_{V \o H^*} = \id$, which implies the assertion by Proposition~\ref{Exponent}.

\subsection[Cauchy's theorem]{} \label{Cauchy}
In 1844, A.-L.~Cauchy proved that a finite group contains an element of order~$p$ for every prime~$p$ that divides the order of the group (cf.~\cite{Cauchy2}, \S~XII, Thm.~5, p.~250). Since the exponent of a finite group is the least common multiple of the orders of all its elements, this can be reformulated by saying that a prime that divides the order of a group also divides its exponent. In this formulation, Cauchy's theorem carries over to semisimple Hopf algebras, as we proceed to prove now. We will need the following lemma:
\begin{lemma}
Suppose that~$p$ is a prime, that $V$\index{$V$} is a finite-dimensional vector space, and that $f: V \rightarrow V$\index{$f$} is an endomorphism of order~$p$.
If $\Tr(f)$ is an integer, we have
$$\Tr(f) \equiv \dim(V) \pmod{p}$$
\end{lemma}
\begin{pf}
In our algebraically closed base field~$K$ of characteristic zero, we choose a primitive $p$-th root of unity~$\zeta$\index{$\zeta$}. The ring of algebraic integers in the cyclotomic field $\Q_p = \Q(\zeta) \subset K$ is then
exactly the ring~$\Z[\zeta]$ (cf.~\cite{Wash}, Prop.~1.2, p.~1), and the principal ideal $P:=(1-\zeta)$ is the unique prime ideal of~$\Z[\zeta]$ lying above $(p) \subset \Z$ (cf.~\cite{Wash}, Lem.~1.4, p.~2). If we denote the dimensions of the eigenspaces of~$f$ by\index{$a_i$}
$$a_i = \dim(\{v \in V \mid f(v) = \zeta^i v\})$$
we find that 
$$\Tr(f) = \sum_{i=0}^{p-1} a_i \zeta^i \qquad \dim(V) = \sum_{i=0}^{p-1} a_i$$
since $f$ is diagonalizable. This shows that $\Tr(f) \equiv \dim(V) \pmod{P}$. Since we have assumed that both numbers are integers, we also get that 
$\Tr(f) \equiv \dim(V)$ modulo $P \cap \Z = (p)$.
\qed
\end{pf}
We note that variants of this lemma, which certainly dates back to the first origins of algebraic number theory, have already earlier been applied to Hopf algebra theory (cf.~\cite{AnSchn3}, Lem.~2.6, p.~433; \cite{SoYp}, Par.~9.6, p.~137). Another such application is the proof of our version of Cauchy's theorem for Hopf algebras, a result that was conjectured by P.~Etingof and S.~Gelaki (cf.~\cite{EG1}, Qu.~5.1, p.~138). We note that the result is known in the case $p=2$ (cf.~\cite{KSZ}, Cor.~4, p.~93).
\begin{thm}
Suppose that $p$ is a prime that divides $\dim(H)$. Then $p$ also divides~$\exp(H)$.
\end{thm}
\begin{pf}
Assume on the contrary that $p$ does not divide $\exp(H)$. From the first statement in Corollary~\ref{SweedlPowInt}, we then get that
$\Lambda^{[p]} = \Lambda^{[1,p]} = \Lambda$. If $\chi_R$ is the character of the regular representation of~$H$, we therefore find that its $p$-th Frobenius-Schur indicator is
$$\nu_p(\chi_R) = \chi_R(\Lambda^{[p]}) = \chi_R(\Lambda) = 1$$
Consider now the endomorphism~$\alpha$ of~$H^{\o p}$ introduced in Paragraph~\ref{FrobSchur}. As we saw there, it preserves the space $(H^{\o p})^H$
of invariants, a space which has, as we discussed in Paragraph~\ref{FrobSchurRegRep}, the dimension~$\dim(H)^{p-1}$. But now Corollary~\ref{FrobSchur}, the first formula for the Frobenius-Schur indicators, together with the preceding lemma yields modulo~$p$ that 
$$1 = \nu_p(\chi_R) = \Tr(\alpha \mid_{(H^{\o p})^H}) \equiv 
 \dim(H)^{p-1} \equiv 0 \pmod{p}$$
This is a contradiction.
\qed
\end{pf}
One could conjecture that Cauchy's theorem carries over to Hopf algebras in an even stronger form, namely as the statement that a semisimple Hopf algebra whose dimension is divisible by a prime~$p$ also contains a grouplike element of order~$p$. However, this is clearly false: As any finite simple group, if it is not abelian, does not have a nontrivial one-dimensional representation, its dual group ring does not contain any nontrivial grouplike elements at all.

Furthermore, one may ask to what extent the above result depends on the assumptions on the base field. We have assumed that it is algebraically closed of characteristic zero, and have used both of these assumptions in the proof. But the assumption that the base field is algebraically closed is clearly not needed, since it is always possible to extend the scalars to the algebraic closure of the base field. However, the result also holds for semisimple and cosemisimple Hopf algebras over fields of positive characteristic; this is a consequence of the lifting theorem of P.~Etingof and S.~Gelaki:
\begin{corollary}
Suppose that $H$ is a semisimple and cosemisimple Hopf algebra over a base field~$K$ of positive characteristic, and that $p$ is a prime that divides $\dim(H)$. Then $p$ also divides~$\exp(H)$.
\end{corollary}
\begin{pf}
As we have just explained, we can assume that the base field is algebraically closed. Since~$K$ is then perfect, there exists a complete discrete valuation ring~$R$ of characteristic zero with residue field~$K$ whose maximal ideal is generated by~$p$, namely the ring of Witt vectors of~$K$ (cf.~\cite{Hasse}, Chap.~10, Sec.~4, p.~156). Such a discrete valuation ring is unique up to isomorphism, and we denote its quotient field by~$F$. By \cite{EG3}, Thm.~2.1, p.~855, there is a free $R$-Hopf algebra~$A$\index{$A$} with $\rank_R(A)= \dim_K(H)$ such that~$A_F:=A \o_R F$\index{$A_F$} is a semisimple and cosemisimple Hopf algebra over the field~$F$ of characteristic zero and that $A/pA$ is isomorphic to~$H$, so that we can actually assume that $H = A/pA$. As we have $\dim_F(A_F) = \dim_K(H)$, our assertion will follow from the preceding theorem if we can show that the exponent of~$A_F$ is equal to the exponent of~$H$. For this, note first that it is obvious that the exponent of~$A_F$ is equal to the exponent of~$A$, and that the exponent of~$H$ divides this number. To see that these numbers are actually equal, we argue as follows: Let $b_1,\ldots,b_n$ be an \mbox{$R$-basis} of the free $R$-module~$A$. If we denote $b_i \o_R 1$ also by~$b_i$, these elements also form an $F$-basis of~$A_F$. The dual basis elements $b_1^*,\ldots,b_n^* \in A_F^*$ map~$A$ to~$R$ and $pA$ to $pR$, and therefore induce mappings 
$\bar{b}_i^* \in \Hom_K(H,K)=H^*$, which are in fact the dual basis elements of the basis of~$H$ consisting of the cosets~$\bar{b}_i$ of~$b_i$ modulo~$pA$.

If $m$ is the exponent of~$A_F$, then it follows from \cite{LR1}, Thm.~2, p.~194 combined with Theorem~\ref{Exponent} that the characteristic of~$K$ does not divide~$m$. Therefore, $K$ contains $m$ distinct $m$-th roots of unity. It is a  consequence of Hensel's lemma (cf.~\cite{Hasse}, Chap.~10, Sec.~7, p.~169) that every $m$-th root of unity can be lifted to~$R$, so that~$F$ also contains $m$ distinct $m$-th roots of unity, all of which are contained in~$R$. Now, since the canonical tensor 
$$t := \sum_{i=1}^n b_i \o b_i^* \in A_F \o A_F^*$$
has order~$m$ by Proposition~\ref{Exponent}, there is a primitive $m$-th root of unity $\zeta \in R$ that is an eigenvalue of the left multiplication by~$t$ in $A_F \o A_F^*$. If we represent the left multiplication by~$t$ by a matrix with respect to the basis $b_i \o b_j^*$ of~$A_F \o A_F^*$, we see that there is a corresponding eigenvector whose components are also in~$F$. Clearing denominators, we can achieve that these components are actually in~$R$, but since~$R$ is a discrete valuation ring with a maximal ideal generated by~$p$, we can assume furthermore that not all components are divisible by~$p$. Then the reduction of this vector modulo~$p$ is still nonzero, which implies that~$\bar{\zeta} \in K$ is still an eigenvalue of the left multiplication by the canonical tensor
$\bar{t} := \sum_{i=1}^n \bar{b}_i \o \bar{b}_i^* \in H \o H^*$. Since both~$F$ and~$K$ contain~$m$ distinct $m$-th roots of unity, 
$\bar{\zeta} \in K$ must still be primitive, which shows that that the order of~$\bar{t}$, which is the exponent of~$H$, is exactly~$m$.
\qed
\end{pf}

\section{The Order} \label{Sec:Order}
\subsection[Order and multiplicity]{} \label{OrderMult}
We now bring in two new concepts, namely the order and the multiplicity of a module:
\begin{defn}
Suppose that~$V$ is an $H$-module. The smallest natural number~$m$ such that $V^{\o m}$ contains a nonzero invariant 
subspace~$(V^{\o m})^H$ is called the order of~$V$ and is denoted by~$\ord(V)$\index{$\ord(V)$}. If~$m$ is the order of~$V$, then the dimension of this invariant subspace is called the multiplicity of~$V$ and is denoted by $\mult(V):=\dim((V^{\o m})^H)$.
\end{defn}
Note that it would a priori be possible that no tensor power of~$V$ contains a nonzero invariant subspace, in which case the order of~$V$ would be infinite. However, we shall prove in Paragraph~\ref{DivThm} that this is not the case.

The notion of the order of a module generalizes the notion of the order of an element in the theory of finite groups. To see this, let~$G$ be a finite group and consider the ring~$K^G$\index{$K^G$} of functions on~$G$, which is isomorphic to the dual group ring~$K[G]^*$. Since this is a commutative Hopf algebra, all its simple modules are one-dimensional, and their characters are given by evaluating a function at a fixed element of the group. This sets up a one-to-one correspondence between the elements of the group and the simple $K^G$-modules under which the product of two group elements corresponds to the tensor product of the modules. Since all these modules are one-dimensional, a tensor power contains a nonzero invariant subspace if and only if it is trivial, which means that the corresponding element of the group is the unit element. Therefore, the order of the module in the sense of the above definition coincides with the order of the element in the sense of group theory. In the case of a general semisimple Hopf algebra~$H$, we see by the same reasoning that, if~$V$ is one-dimensional, and therefore determined by its character $\gamma: H \rightarrow K$, the order of~$V$ in the sense of the above definition coincides with the order of~$\gamma$ in the group~$G(H^*)$ of grouplike elements of the dual Hopf algebra~$H^*$.

\subsection[The divisibility theorem]{} \label{DivThm}
To proceed further, we will need some properties of symmetric polynomials.
Consider the polynomial ring~$\Q[x_1,x_2,\ldots,x_n]$ in $n$ variables. Recall that, for a nonnegative integer~$k$, the $k$-th elementary symmetric polynomial $e_k=e_k(x_1,\ldots,x_n)$ is defined as
$$e_k(x_1,\ldots,x_n) := 
\sum_{i_1<i_2<\ldots<i_k} x_{i_1} x_{i_2} \ldots x_{i_k}$$
In particular, we have
$e_0=1$, $e_1=x_1 + x_2 + \ldots + x_n$, $e_n = x_1 x_2 \ldots x_n$, and $e_k=0$ for $k>n$. The fundamental theorem on symmetric polynomials (cf.~\cite{SchejaStorch2}, Satz~54.13, p.~57) asserts that every symmetric polynomial with integer coefficients can be expressed as a polynomial with integer coefficients in the elementary symmetric polynomials. This holds in particular for the power sums
$$s_k = s_k(x_1,\ldots,x_n) := \sum_{i=1}^n x_i^k$$
Conversely, the power sums have the property that every symmetric polynomial can be expressed as a polynomial with rational coefficients in the power sums. The polynomials that achieve this for the elementary symmetric polynomials are the so-called (fractional) Newton polynomials~$Q_n$\index{$Q_n$}, which are defined via the following $n \times n$-determinant:
$$
Q_n(x_1,\ldots,x_n) := \frac{(-1)^n}{n!}
\begin{vmatrix}
x_1 & 1 & 0 &  0 & \ldots & 0 \\
x_2 & x_1 & 2 & 0 & \ldots & 0 \\
x_3 & x_2 & x_1 & 3 & \ldots & 0 \\
\vdots & \vdots & \vdots & \vdots &  & \vdots \\
x_{n-1} & x_{n-2} & x_{n-3} & x_{n-4 }&  \ldots & n -1 \\
x_n & x_{n-1} & x_{n-2} & x_{n-3 }&  \ldots & x_1
\end{vmatrix}
$$
The polynomials $Q_1,\ldots,Q_{n-1}$, which involve only fewer variables, can of course also be considered as elements of~$\Q[x_1,\ldots,x_n]$.
The formula that expresses the elementary symmetric polynomials in terms of the power sums is known as Newton's formula (cf.~\cite{SchejaStorch2}, Satz~57.9, p.~110):
\begin{lemma}
For $n \ge 1$ and $j=1,\ldots,n$, we have
$e_j = (-1)^j Q_j(s_1,s_2,\ldots, s_n)$.
\end{lemma}
With the help of this formula, we can now prove the following theorem:
\begin{thm}
The order of a nonzero $H$-module $V$ is a finite number. It is not larger than $\dim(H)$ and divides $\dim(H) \mult(V)$.
\end{thm}
\begin{pf}
Recall from Paragraph~\ref{SecondForm} that, for any $\varphi \in H^*$, the trace of the left multiplication 
$$L_\varphi: H^* \rightarrow H^*,~\psi \mapsto \varphi \psi$$
is $n \varphi(\Lambda)$, where $n= \dim(H)$. Let $\chi \in H^*$ be the character of~$V$, and let 
$\lambda_1, \ldots, \lambda_n$ be the not necessarily distinct eigenvalues of~$L_\chi$. If $m:=\ord(V) \le n$, we have for $k=1,\ldots,m-1$ that 
$$s_k(\lambda_1,\ldots,\lambda_n) = \Tr(L_{\chi^k}) = n \chi^k(\Lambda) = 0$$
since $V^{\o k}$ does not contain a nonzero invariant subspace. For $k=m$, we get similarly that $s_m(\lambda_1,\ldots,\lambda_n) = n \mult(V)$. By Newton's formula, we therefore have 
$$
e_m(\lambda_1,\ldots,\lambda_n) = 
\frac{1}{m!}
\begin{vmatrix}
0 & 1 & 0 &  0 & \ldots & 0 \\
0 & 0 & 2 & 0 & \ldots & 0 \\
0 &  & 0 & 3 & \ldots & 0 \\
\vdots & \vdots & \vdots & \vdots &  & \vdots \\
0 & 0 & 0 & 0 &  \ldots & m-1 \\
s_m & 0 & 0 & 0&  \ldots & 0
\end{vmatrix}
= \frac{(-1)^{m-1}}{m} s_m(\lambda_1,\ldots,\lambda_n)$$
Now note that the eigenvalues $\lambda_1,\ldots,\lambda_n$, and therefore the elementary symmetric functions~$e_k(\lambda_1,\ldots,\lambda_n)$, are algebraic integers. To see this, consider the left multiplication by~$\chi$ not on the whole dual Hopf algebra~$H^*$, but only on the 
character ring~$\Ch(H)$\index{$\Ch(H)$}. Its matrix representation with respect to the basis consisting of the irreducible characters has integer entries, and therefore the Cayley-Hamilton theorem implies that it satisfies a monic polynomial with integer coefficients, namely its characteristic polynomial. Evaluating this on the unit of the character ring, we see that $\chi$ itself satisfies this polynomial, and therefore also~$L_\chi$ satisfies 
this polynomial. Since $\lambda_1,\ldots,\lambda_n$ are roots of this polynomial, they must be algebraic integers. This shows that the fraction
$$\frac{n \mult(V)}{m} = \frac{1}{m} s_m(\lambda_1,\ldots,\lambda_n)=
(-1)^{m-1}e_m(\lambda_1,\ldots,\lambda_n)$$
is an algebraic integer, which therefore must be an integer (cf.~\cite{SchejaStorch2}, \S~56, Bsp.~3, p.~91). This proves the divisibility assertion.

It still remains to be shown that the order is finite and bounded by~$n$. So suppose that this is not the case. The reasoning above then shows that the power sums $s_1,\ldots,s_n$ of the eigenvalues are zero, and therefore, by Newton's formula, the elementary symmetric functions~$e_1,\ldots,e_n$ of the eigenvalues are zero. From this, we conclude as follows that all eigenvalues $\lambda_1,\ldots,\lambda_n$ are zero: First, we have 
$e_n(\lambda_1,\ldots,\lambda_n)=\lambda_1 \lambda_2 \ldots \lambda_n=0$, which implies that~$\lambda_i=0$ for some~$i$. Considering the next elementary symmetric function, we have
$$e_{n-1}(\lambda_1,\ldots,\lambda_n) = \lambda_1 \cdot \ldots \cdot\lambda_{i-1} \lambda_{i+1} \cdot \ldots \cdot \lambda_n=0$$
since all the other summands in the definition of~$e_{n-1}$ vanish. This implies that we have $\lambda_j=0$ for another~$j \neq i$. Proceeding in this way, we arrive at the assertion that all eigenvalues of~$L_\chi$ are zero. But this is not the case, as a nonzero integral of~$H^*$ is an eigenvector for~$L_\chi$ corresponding to the eigenvalue~$\dim(V)$.~\qed 
\end{pf}
It should be noted that this reasoning shows that the elementary symmetric functions of the eigenvalues are not only algebraic integers, but in fact  integers. This holds since they can be calculated, with the help of Newton's formula, by rational operations from the power sums. But the power sums are integers,
as they are multiplicities of the trivial module in some tensor power times
the dimension of~$H$. Since the elementary symmetric functions of the eigenvalues are, up to a sign, the coefficients of the characteristic polynomial, we get the following consequence:
\begin{corollary}
The characteristic polynomial of the left multiplication by a character on~$H^*$ has integral coefficients.
\end{corollary}
We note that there exists at least one other proof of this corollary, which is also rather instructive. However, we do not give it here because it would take
us too far afield.

\subsection[An example]{} \label{Examplepq}
We will see in Paragraph~\ref{IndecompMat} that the order of a module is not only bounded by the dimension of the Hopf algebra, but even by the dimension of the character ring. Another conjecturable strengthening of the divisibility theorem above, namely the statement that $\ord(V)$ always divides~$\dim(H)$, is, however, not true, as the following example shows:
Suppose that~$p$\index{$p$} and~$q$\index{$q$} are two prime numbers with the property that~$p$
divides~$q-1$. Choose an element~$a$\index{$a$} of order~$p$ in the group of
units~$\Z_q^\times$, and consider the semidirect 
product~$G := \Z_q \rtimes \Z_p$ defined by 
$$(m,n)(m',n') := (m + a^n m', n + n')$$
For a $q$-th root of unity~$\zeta$\index{$\zeta$}, we can consider the
base field~$K$ as a $\Z_q$-module via~$\zeta$. If $H:=K[G]$ is the group ring, we can form the induced module
$$V:=H \o_{K[\Z_q]} K$$
The elements $v_1,\ldots,v_p$, where~$v_n := (0,n) \o_{K[\Z_q]} 1_K$\index{$v_n$}, then form a basis of~$V$. Reading these indices modulo~$p$, we have $(0,1).v_n=v_{n+1}$ and 
$$(1,0).v_n = (0,n)(a^{-n},0) \o_{K[\Z_q]} 1_K = \zeta^{a^{-n}} v_n$$
Now suppose that~$\zeta$ is primitive. Since the elements 
$a,a^2,\ldots,a^p \in \Z_q^\times$ are all distinct, 
we see from the standard isomorphism $\Z_q^\times \cong \Gal(\Q_q/\Q)$ that also the
elements $\zeta^{a^n}$, for $n \in \Z_p$, are distinct, and therefore $(1,0)$ operates as a diagonalizable operator with distinct eigenvalues. In particular, this shows that~$V$ is simple, since every nonzero $H$-submodule
must on the one hand be stable under this operator, and therefore be spanned by some of the vectors~$v_n$, and must on the other hand be stable under the action of~$(0,1)$, and therefore contain all the vectors~$v_n$ if it contains one of them.

The order of~$V$ can be determined from the following result:
\begin{prop}
The space of invariants~$(V^{\o m})^H$ is nonzero if and only if there are
numbers~$n_1,\ldots,n_m \in \Z_p$ such that 
$a^{n_1} + a^{n_2} + \ldots + a^{n_m} = 0 \in \Z_q$.
\end{prop}
\begin{pf}
An arbitrary tensor $w \in V^{\o m}$ can be written uniquely in the form
$w = \sum_{i_1,\ldots,i_m=1}^p \alpha_{i_1,\ldots,i_m} v_{i_1} \o \ldots \o v_{i_m}$.
We have
\begin{align*}
(1,0).w &=
\sum_{i_1,\ldots,i_m=1}^p \alpha_{i_1,\ldots,i_m}
\zeta^{a^{-i_1} + a^{-i_2} + \ldots + a^{-i_m}} v_{i_1} \o \ldots \o v_{i_m} \\
(0,1).w &=
\sum_{i_1,\ldots,i_m=1}^p \alpha_{i_1,\ldots,i_m} v_{i_1+1} \o \ldots \o v_{i_m+1}
\end{align*}
Therefore, if $w$ is a nonzero invariant tensor, then we have $\alpha_{i_1,\ldots,i_m} \neq 0$ for
some elements $i_1,\ldots,i_m \in \Z_p$, which implies that
$\zeta^{a^{-i_1} + a^{-i_2} + \ldots + a^{-i_m}} = 1$ or
$$a^{-i_1} + a^{-i_2} + \ldots + a^{-i_m} = 0 \in \Z_q$$
Conversely, suppose that
$a^{n_1} + a^{n_2} + \ldots + a^{n_m} = 0 \in \Z_q$. Then we also have
$a^{n_1-i} + a^{n_2-i} + \ldots + a^{n_m-i} = 0 \in \Z_q$, and therefore the tensor
$$w:=\sum_{i=1}^p v_{i-n_1} \o \ldots \o v_{i-n_m}$$
is invariant under the action of~$(1,0)$ and~$(0,1)$, and therefore invariant under the action
of~$H$. 
\qed
\end{pf}
Let us consider a concrete case: For~$p=5$, $q=11$, we can choose $a=3$. The powers of~$3$ modulo~$11$ are
$1,3,9,5$, and~$4$. Since no sum of two of these is zero modulo~$11$, we see that $\ord(V)>2$.
On the other hand, we have $1+5+5 \equiv 0 \pmod{11}$. Therefore, we have $\ord(V)=3$, which does not divide $\dim(H)=55$.

\subsection[The dimension of the simple modules]{} \label{DimSimMod}
A special case of the divisibility theorem given in Paragraph~\ref{DivThm} 
is the following result (cf.~\cite{KSZ}, Thm.~4, p.~91):
\begin{corollary} 
If $H$ has a non-trivial self-dual simple module, then the dimension of~$H$ is even.
\end{corollary}
This holds because, as we discussed in Paragraph~\ref{FrobSchurThm}, a simple module~$V$ is self-dual if and only if $(V \o V)^H \neq 0$, and therefore a non-trivial self-dual simple module has order~$2$ and multiplicity~$1$. Note that this proof, in contrast to the one given in~\cite{KSZ}, does not rely on the notion of the exponent. Let us give a third proof, which is also very short and uses the exponent instead of the divisibility theorem: If $\chi$ is the character of a nontrivial self-dual simple module, we have $\chi(\Lambda^{[2]})= \pm 1$ by the Frobenius-Schur theorem~\ref{FrobSchurThm} and $\chi(\Lambda)= 0$ since the module is nontrivial. But if the dimension of~$H$ were odd, the exponent of~$H$ would be odd by Theorem~\ref{Exponent}, and therefore we would have $\Lambda = \Lambda^{[1,2]} = \Lambda^{[2]}$ by~Corollary~\ref{SweedlPowInt}, a contradiction.

The remarkable fact that it is possible to give a proof of this result by two apparently different techniques may be partially explained by the circumstance that the Frobenius-Schur theorem implies that for an irreducible character we have $\chi^{[2]}(\Lambda)= \pm \chi^2(\Lambda)$. It seems that a corresponding relation between the Sweedler powers~$\chi^{[m]}$ and the ordinary powers~$\chi^m$ of a character~$\chi$ is not known at present.

We have already explained in~\cite{KSZ}, Thm.~5, p.~93 that the above corollary implies that the dimension of~$H$ must be even if~$H$ has a simple module of even dimension. Here, we give another application:
\begin{thm}
Suppose that the dimension of~$H$ is odd. If $H$ has a simple module of dimension~$3$, then the dimension of~$H$ is divisible by~$3$.
\end{thm} 
\begin{pf}
Suppose that~$V$ is a three-dimensional simple module with character~$\chi$. By Schur's lemma, the trivial module appears in the decomposition 
of~$V^* \o V$. As we have just said, the simple constituents of~$V^* \o V$
cannot have even dimension, and therefore their possible dimensions are~$1$, $3$, $5$, and~$7$. If~$V^* \o V$ had a unique simple constituent of dimension~$7$, this constituent would be self-dual, which is impossible by the corollary above. But for dimension considerations, there cannot be two simple constituents of dimension~$7$, and therefore simple constituents of dimension~$7$ cannot appear at all. For the same reason, simple constituents of dimension~$5$ cannot appear, and the number of three-dimensional simple constituents is either~$0$ or~$2$.

Therefore, $V^* \o V$ decomposes either into nine one-dimensional simple modules or into three one-dimensional and two three-dimensional simple modules. In both cases, the number of one-dimensional simple constituents is divisible by~$3$. Now a one-dimensional module with character~$\gamma$ appears in~$V^* \o V$ if and only if $\chi=\chi\gamma$, and if it appears, it appears with multiplicity one (cf.~\cite{NR}, Thm.~10, p.~303). The one-dimensional constituents of~$V^* \o V$ therefore form a subgroup of~$G(H^*)$, namely the isotropy group of~$\chi$ under right multiplication by elements of~$G(H^*)$, whose order is divisible by~$3$. By the Nichols-Zoeller theorem (cf.~\cite{M}, Thm.~3.1.5, p.~30), this implies that the dimension of~$H$ is divisible by~$3$.
\qed
\end{pf}
The above theorem is a variant of a result of S.~Burciu (cf.~\cite{Burciu}, Cor.~8, p.~93). His setup deviates slightly from ours: He does not assume that the characteristic of the base field is zero, but adds the assumption that~$H$ has no even-dimensional simple modules, which is, as pointed out above, automatic in the characteristic zero case. His proof is rather different and mimics the methods from~\cite{NR}.

\section{The Index} \label{Sec:Index}
\subsection[Indecomposable matrices]{} \label{IndecompMat}
Suppose that $\chi_1,\ldots,\chi_k$\index{$\chi_i$} are the distinct irreducible characters of the semisimple Hopf algebra~$H$ under consideration, and that $V_1,\ldots,V_k$\index{$V_i$} are simple $H$-modules of dimension~$n_1,\ldots,n_k$\index{$n_i$} corresponding to these characters. We can assume that~$V_1=K$, the base field considered as a trivial $H$-module, with character $\chi_1=\varepsilon$, the counit. If $V$ is an arbitrary $H$-module with character~$\chi$, we have already used in the proof of Theorem~\ref{DivThm} that the matrix representation of the left multiplication by~$\chi$ on the character ring~$\Ch(H)$ with respect to the basis $\chi_1,\ldots,\chi_k$ has nonnegative integer entries: If
$$\chi \chi_j = \sum_{i=1}^k a_{ij} \chi_i$$
then $a_{ij}$ is just the multiplicity of~$V_i$ in the decomposition 
of~$V \o V_j$ into simple modules---clearly a nonnegative integer. We are therefore in the position to apply the theory of nonnegative matrices, which we do in this section.

Recall that the $k \times k$-matrix $A=(a_{ij})_{i,j=1,\ldots,k}$\index{$A$} is called decomposable if it is possible to find a decomposition $I_k = M \cup N$ of  $I_k=\{1,\ldots,k\}$ into disjoint nonempty sets~$M$\index{$M$} and~$N$\index{$N$} such that $a_{ij}=0$ whenever $i \in M$ and $j \in N$ (cf.~\cite{Gant}, Abschn.~13.1., Def.~2, p.~395); otherwise, it is called indecomposable. By changing the enumeration of the indices in such a way that the indices in~$M$ come first and the indices in~$N$ come second, we can achieve that a decomposable matrix has zeros in the upper right rectangular part. Note that any power $A^m$ of a decomposable matrix is still decomposable, since its matrix elements are
$$\sum_{i_1,\ldots,i_{m-1}=1}^k a_{ii_1} a_{i_1 i_2} a_{i_2 i_3} \ldots
a_{i_{m-2} i_{m-1}} a_{i_{m-1}j}$$
and these are zero whenever $i \in M$ and~$j \in N$.

To link this notion to Hopf algebras, we need another ingredient. Suppose that $V$ is an $H$-module. The set of all elements of~$H$ that act on~$V$ identically as zero is called the annihilator of~$V$. This is obviously a two-sided ideal, but not necessarily a Hopf ideal. But if we denote by~$J$\index{$J$} the intersection of the annihilators of all the tensor powers~$V^{\o m}$ of~$V$, including the trivial module~$K$ for $m=0$, we get the largest Hopf ideal contained in the annihilator:
\begin{lemma}
$J$ is a Hopf ideal of~$H$. Every other Hopf ideal of~$H$ that is contained in the annihilator of~$V$ is contained in~$J$.
\end{lemma}
For a proof of this lemma, we refer to~\cite{Rie}, Thm.~1, p.~125 (see also~\cite{PQ}, Prop.~1*, p.~328). Note that the usage of the notion of a Hopf algebra in~\cite{Rie} differs from our usage. To adopt the proof, one has to use the fact that in a finite-dimensional Hopf algebra a bi-ideal is a Hopf ideal (cf.~\cite{NichQuot}, Thm.~1, p.~1791; see also~\cite{PQ}, Lem.~6, p.~331). Since a simple $H$-module can be embedded into a tensor power of~$V$ if and only if it is annihilated by~$J$, this also shows that if a simple module can be embedded into a tensor power of~$V$, its dual can also be embedded into a tensor power of~$V$.

Now, let as before~$\chi$ be the character and $I$ be the annihilator of~$V$, and denote by~$A$ the matrix representation of the left multiplication by~$\chi$ on the character ring~$\Ch(H)$ with respect to the basis consisting of the irreducible characters~$\chi_1,\ldots,\chi_k$. The following result connects these things:
\begin{prop}
The following statements are equivalent:
\begin{enumerate}
\item $I$ does not contain a nonzero Hopf ideal.
\item $A$ is indecomposable.
\end{enumerate}
\end{prop}
\begin{pf}
The statement is correct if~$V$ is zero, even if~$H$ is one-dimensional, so let us assume that~$V$ is nonzero. We first show that the first statement implies the second. Let us assume that~$A$ is decomposable, and choose a corresponding decomposition $I_k = M \cup N$. Since we saw above that powers of a decomposable matrix are still decomposable, we then have for $i \in M$ and~$j \in N$ that~$V_i$ never appears as a direct summand of~$V^{\o m} \o V_j$. On the other hand, it follows from the above lemma that the intersection of the annihilators of the tensor powers of~$V$ are zero, which means that every simple module appears as a constituent of some~$V^{\o m}$. Suppose now that $V_i$ appears as a constituent 
of~$V^{\o r}$ and that $V_j^*$ appears as a constituent of~$V^{\o s}$, which means that there are injective $H$-linear maps from $V_i$ to~$V^{\o r}$ and from $V_j^*$ to~$V^{\o s}$, which we indicate by hooked arrows. By Schur's lemma, the trivial module~$K$ appears as a
constituent of $V_j^* \o V_j$:
$$K \hookrightarrow V_j^* \o V_j  \hookrightarrow V^{\o s} \o V_j$$
This implies that $V_i$ appears as a constituent of $V^{\o (r+s)} \o V_j$:
\begin{samepage}
$$V_i \hookrightarrow V^{\o r} \o K \hookrightarrow 
V^{\o r} \o V^{\o s} \o V_j = V^{\o (r+s)} \o V_j$$
We have therefore reached a contradiction.
\end{samepage}

Let us next prove that the second statement implies the first; so assume that~$A$ is indecomposable. Define the sets
\begin{align*}
&M:=\{i \le k \mid V_i \; \text{cannot be embedded into} \; V^{\o m} \; \text{for any} \; m\} \\
&N:=\{j \le k \mid V_j \; \text{can be embedded into} \; V^{\o m} \; \text{for some} \; m\}
\end{align*}
If $j \in N$, choose $m$ such that
$V_j$ can be embedded into~$V^{\o m}$. $V_i$ then appears in $V \o V_j$ with multiplicity~$a_{ij}$,
and therefore $V_i$ appears in $V^{\o (m+1)}$ at least with multiplicity~$a_{ij}$. Therefore, if~$i \in M$, we have $a_{ij}=0$. Since~$A$ is indecomposable and~$N$ is not empty, as it contains the indices corresponding to the simple constituents of the nonzero module~$V$, $M$ must be empty. This means that every simple module appears as a constituent of some tensor power of~$V$, and therefore the intersection of the annihilators of all these tensor powers must be zero. But by the preceding lemma, this is the largest Hopf ideal contained in~$I$, and the assertion follows.
\qed
\end{pf}
As promised in Paragraph~\ref{Examplepq}, we can now prove that the order of a module is bounded by the dimension of the character ring:
\begin{corollary}
For a nonzero $H$-module $V$, we have $\ord(V) \le \dim(\Ch(H))$.
\end{corollary}
\begin{pf}
As above, let $J$ be the intersection of the annihilators of all the tensor powers of~$V$. Since this is a Hopf ideal, $H/J$ is a Hopf algebra, and~$V$ can also be considered as a module over this algebra. If $\chi'$\index{$\chi'$} denotes the character of~$V$ as a module over~$H/J$, then the matrix~$A$ that represents the left multiplication by~$\chi'$ with respect to the basis consisting of the irreducible characters of~$H/J$ is, by the preceding proposition, indecomposable. If $l:=\dim(\Ch(H/J)) \le \dim(\Ch(H))$, then one of the powers $A, A^2, A^3, \ldots, A^l$ has a nonzero
$(1,1)$-component (cf.~\cite{Gant}, Abschn.~13.1, Folgerung, p.~397). If $A^m$ is this matrix, this means that \mbox{$V^{\o m} \o K$} contains  the trivial module with nonzero multiplicity, i.e., $V^{\o m}$ contains a nonzero invariant submodule.
\qed
\end{pf}
Although this result is not stated as such, the above corollary should be viewed as a variant of Theorem~2 in~\cite{Rie}, p.~126. The proof given there has the advantage not to rely on the theory of nonnegative matrices. Nonetheless, we have given the above proof because it fits well to the material that we will discuss below. Moreover, as pointed out in~\cite{Gant}, loc.~cit., the above argument shows the stronger statement that $\ord(V)$ is in fact bounded by the degree of the minimum polynomial of~$\chi'$.

\subsection[The normal form]{} \label{NormForm}
Even if it is not indecomposable, a matrix with nonnegative real entries can be cast into a certain normal form by changing the enumeration of rows and columns
(cf.~\cite{Gant}, Abschn.~13.4, p.~417). In this normal form, the matrix has a lower block triangular form with indecomposable diagonal blocks. But in our situation, where the matrix~$A$ arises from the left multiplication by a character~$\chi$ of an $H$-module~$V$, we can, as P.~Etingof pointed out (cf.~\cite{EtingPriv}), say even more:
\begin{prop}
If the irreducible characters $\chi_1,\ldots,\chi_k$\index{$\chi_i$} are suitably enumerated, the matrix~$A$ is block diagonal with indecomposable diagonal blocks.
\end{prop}
\begin{pf}
The assertion holds if~$A$ is indecomposable, so suppose now that we have a decomposition $I_k = M \cup N$ of  $I_k=\{1,\ldots,k\}$ into disjoint nonempty sets~$M$\index{$M$} and~$N$\index{$N$} such that $a_{ij}=0$ whenever $i \in M$ and $j \in N$. As we explained in Paragraph~\ref{IndecompMat}, every power of~$A$ also has this property, which implies that $V^{\o m} \o V_j$ does not contain~$V_i$ as a constituent. Now, if $W$ is a simple constituent of~$V$, we saw in the discussion after Lemma~\ref{IndecompMat} that~$W^*$ can be embedded into a tensor power~$V^{\o m}$ of~$V$. Therefore, also $W^* \o V_j$ does not contain~$V_i$ as a constituent. But this implies that $W \o V_i$ does not contain~$V_j$ as a constituent (\cite{NR}, Thm.~9, p.~303; \cite{So4}, Par.~3.6, p.~212). Applying this to all simple constituents of~$V$, we see that $V \o V_i$ does not contain~$V_j$, so that $a_{ji}=0$.

Now consider the submatrix $(a_{ij})_{i,j \in M}$ formed from the indices that belong to~$M$, and also the analogous matrix formed from the indices that belong to~$N$. If these matrices are decomposable, we choose a corresponding decomposition of~$M$, resp.~$N$, into disjoint nonempty subsets. Proceeding in this way, we arrive at a decomposition $I_k = M_1 \cup M_2 \cup \ldots \cup M_s$ into disjoint nonempty sets~$M_r$\index{$M_r$}\index{$s$} with the following two properties: First, the submatrices $(a_{ij})_{i,j \in M_r}$  are indecomposable. Second, for indices $r$ and $t$ satisfying $1 \le r < t \le s$, we have that $a_{ij}=0$ whenever $i \in M_r$ and $j \in M_t$. By applying the preceding argument with $M_1$ and $I_k \setminus M_1$ instead of~$M$ and~$N$, we see  
that $a_{ij}=0$ also if $j \in M_1$ and $i \notin M_1$. This shows that 
$I_k = M_2 \cup I_k \setminus M_2$ is also a decomposition to which the preceding argument applies, so that we get that $a_{ij}=0$ if $j \in M_2$ and 
$i \notin M_2$. Continuing in this way, we get that $a_{ij}=0$ whenever~$i$ and~$j$ belong to distinct sets of the decomposition 
$I_k = M_1 \cup M_2 \cup \ldots \cup M_s$.

If we now enumerate the irreducible characters~$\chi_1,\ldots,\chi_k$ so that the elements in~$M_1$ come first, the elements of~$M_2$ come second, and so on, then the preceding argument shows that the matrix~$A$ corresponding to this enumeration is block diagonal with indecomposable blocks. Note that this proof shows that the enumeration can be changed in such a way that $\chi_1$ is still the counit.
\qed
\end{pf}

\subsection[The Perron-Frobenius theorem]{} \label{Index}
For the following, we need to recall the main result of the theory of nonnegative matrices, namely the Perron-Frobenius theorem (cf.~\cite{Gant}, Abschn.~13.2, Satz~2, p.~398). We divide the theorem into two parts, the first part saying the following:
\begin{thm}[Part 1]
Suppose that $A$ is an indecomposable square matrix with nonnegative real entries. Then $A$ has a positive eigenvalue~$\lambda$, called the Perron-Frobenius eigenvalue, with the property that $|\mu| \le \lambda$ for every other eigenvalue~$\mu$. The algebraic multiplicity of~$\lambda$ is one, i.e., $\lambda$ is a simple root of the characteristic polynomial. The corresponding eigenvector, which is therefore unique up to scalar multiples, can be chosen to have positive components. Such an eigenvector is then called a Perron-Frobenius eigenvector.
\end{thm}
For the second part, we deviate a little bit from the formulation in~\cite{Gant}, although everything we state is proved there.
We define the index of imprimitivity, or briefly the index of~$A$\index{$\ind(A)$}, to be
$$\ind(A) := 
|\{\mu \mid \mu~\text{is an eigenvalue of}~A~\text{with}~|\mu|=\lambda \}|$$
the number of eigenvalues for which the above inequality is actually an equality (cf.~\cite{Gant}, Abschn.~13.5, Def.~3, p.~422). If $\zeta \in \C$ is a primitive $\ind(A)$-th root of unity, then the second part of the Perron-Frobenius theorem can be formulated in the following way (cf.~\cite{Gant}, Abschn.~13.2, Gl.~(31), p.~404):
\begin{thm}[Part 2]
There is a diagonal matrix~$D$\index{$D$} whose diagonal entries are $\ind(A)$-th roots of unity such that $D A D^{-1} = \zeta A$.
\end{thm}
In particular, since $A$ is similar to~$\zeta A$, $\zeta \mu$ is an eigenvalue of~$A$ whenever~$\mu$ is, showing that the eigenvalues~$\mu$ of~$A$ that satisfy~$|\mu|=\lambda$ are exactly the numbers of the form~$\mu=\zeta^m \lambda$. Furthermore, if $x$ is an eigenvector of~$D$ corresponding to the eigenvalue $\zeta^m$, then $Ax$ is an eigenvector of~$D$ corresponding to the eigenvalue~$\zeta^{m+1}$, so that~$A$ shifts the eigenspaces of~$D$ around cyclicly, although, since it is not necessarily invertible, it does not always induce an isomorphism between these eigenspaces.

We now apply this theorem to Hopf algebras. For this, let us briefly explain how to deal with the real numbers that appear as the components of a Perron-Frobenius eigenvector and are in general not elements of our base field~$K$. Note first that the $\Q$-vector space $\Ch_\Q(H):=\Span_\Q(\chi_1,\ldots,\chi_k)$\index{$\Ch_\Q(H)$} is closed under multiplication, because, as we saw above, the structure constants for the products of the irreducible characters are nonnegative integers. From this \mbox{$\Q$-algebra}, we can extend the scalars to the real numbers and work in the algebra $\Ch_\R(H) := \R \o_\Q \Ch_\Q(H)$\index{$\Ch_\R(H)$}, in which linear combinations of characters with real coefficients make sense.

Suppose now that~$V$ is an $H$-module with character~$\chi$. We have already seen in Lemma~\ref{IndecompMat} that the annihilator of~$V$ contains a unique largest Hopf ideal~$J$, namely the intersection of the annihilators of the tensor powers of~$V$. As in the proof of Corollary~\ref{IndecompMat}, we consider~$V$ as a module over the quotient Hopf algebra~$H/J$, and it follows from Proposition~\ref{IndecompMat} that the matrix representation~$A$ of the left multiplication by the character of~$V$ on the character ring~$\Ch(H/J)$ with respect to the basis that consists of the irreducible characters of~$H/J$ is indecomposable. Our first task is to determine its Perron-Frobenius eigenvalue:
\begin{prop}
The Perron-Frobenius eigenvalue of~$A$ is~$\dim(V)$.
\end{prop}
\begin{pf}
The following argument is taken from~\cite{So4}, Thm.~3.7, p.~213; we repeat it here for the sake of completeness. We first note that, by replacing~$H$ by~$H/J$, we can assume that $J$ is zero. The matrix $A=(a_{ij})_{i,j=1,\ldots,k}$ is then defined by the equation $\chi \chi_j = \sum_{i=1}^k a_{ij} \chi_i$. Now 
let $\lambda>0$ be the Perron-Frobenius eigenvalue of~$A$, and choose a corresponding Perron-Frobenius eigenvector $(x_1,\ldots,x_k)$ with positive real entries. Working in the algebra~$\Ch_\R(H)$, we then have the equation 
$\chi \chi' = \lambda \chi'$, where 
$\chi' := \sum_{i=1}^k x_i \chi_i$\index{$\chi'$}.
Evaluating this equation at the unit element, we get
$$\dim(V) \sum_{i=1}^k x_i n_i = \chi(1) \chi'(1) 
= \lambda \sum_{i=1}^k x_i n_i$$
Since the $x_i$ are positive, the term $\sum_{i=1}^k x_i n_i$ is nonzero, and therefore we get $\lambda = \dim(V)$.
\qed
\end{pf}
It may be noted that the largest eigenvalue of the left multiplication by a character is always its degree, whether the associated matrix is indecomposable or not. This holds since a corollary of the Perron-Frobenius theorem (cf.~\cite{Gant}, Abschn.~13.3, Satz~3, p.~409) still yields for any matrix with nonnegative real entries a nonnegative eigenvalue that is an upper bound for the absolute values of all the other eigenvalues. Furthermore, there is a corresponding eigenvector with nonnegative real entries, which is the only thing that was needed in the above argument. In our situation, it is easy to find such an eigenvalue and a corresponding eigenvector: The character of the regular representation is an integral (cf.~\cite{LR2}, Prop.~2.4, p.~273), and multiplying this integral by the given character shows that the corresponding eigenvalue of the character is its degree. As the components of this integral in the basis $\chi_1,\ldots,\chi_k$ are exactly the dimensions $n_1,\ldots,n_k$, which are all strictly positive, we in fact see that the indecomposable diagonal blocks that we discussed in Paragraph~\ref{NormForm} all have the Perron-Frobenius eigenvalue~$\dim(V)$ (cf.~\cite{Gant}, Abschn.~13.4, Satz~7, p.~419).

\subsection[The index formula]{} \label{IndForm}
After having determined the first invariant, we turn to the second invariant.
Note first that~$V$ is a faithful $G(H/J)$-module: If $g \in G(H/J)$ is a grouplike element of the quotient Hopf algebra~$H/J$ that acts as the identity on~$V$, then it also acts as the identity on every tensor power of~$V$, which implies that~$g-1 \in J/J$, so that $g=1$. We now consider the subgroup of~$G(H/J)$ consisting of those elements that act on~$V$ as a scalar multiple of the identity:
\begin{defn}
$G_V:=\{g \in G(H/J) \mid \exists \, \xi \in K \; \forall \, v \in V: 
g.v = \xi v\}$\index{$G_V$}\index{$\xi$}
\end{defn}
Using this group, we can give the following formula for the index, which generalizes a well-known result in the representation theory of finite groups (cf.~\cite{Isaacs}, Lem.~(2.27), p.~27):
\begin{thm}
The group $G_V$ is cyclic and contained in the center~$Z(H/J)$\index{$Z(H)$} of~$H/J$. Its order is equal to the index of imprimitivity of~$A$\index{$\ind(A)$}: 
$$\ind(A) = |G_V|$$
\end{thm}
\begin{pf}
As before, by replacing $H$ by~$H/J$, we can assume that $J$ is zero. 
For~$g \in G_V$, there is by definition a unique scalar~$\xi \in K$ such 
that~$g.v = \xi v$ for every~$v \in V$. In this way, we get a group homomorphism
$$G_V \rightarrow K^\times,~g \mapsto \xi$$
into the group of nonzero elements~$K^\times$ of~$K$, which, as we saw above, is injective. This shows that~$G_V$ is isomorphic to a finite subgroup of~$K^\times$ and therefore cyclic. To see that it is central, note that, for~$g \in G_V$ and~$h \in H$, the elements~$gh$ and~$hg$ act in the same way on every tensor power of~$V$, so that~$gh-hg \in J$ and therefore~$gh=hg$. 

It remains to show the index formula. Using a bar to denote the character of the dual module, we have the identity
$$\sum_{i=1}^k \chi \chi_i \o \bar{\chi}_i = 
\sum_{i=1}^k \chi_i \o \bar{\chi}_i \chi$$
(cf.~\cite{So4}, Prop.~3.5, p.~211). Evaluating the second tensorand 
on an element $g \in G_V$ that acts on~$V$ by multiplication with~$\xi \in K$, we get
$$\chi \sum_{i=1}^k \chi_i(g^{-1}) \chi_i = 
\xi \dim(V)\sum_{i=1}^k \chi_i(g^{-1}) \chi_i$$
which means that $\sum_{i=1}^k \chi_i(g^{-1}) \chi_i$ is an eigenvector for the left multiplication by~$\chi$ corresponding to the 
eigenvalue~$\xi \dim(V)$. Obviously, the absolute value of this eigenvalue is~$\dim(V)$, so that, if $\zeta$ is the primitive $\ind(A)$-th root of unity\index{$\ind(A)$} appearing in the second part of the Perron-Frobenius theorem, the discussion there shows that~$\xi$ is a power of~$\zeta$. Since~$\xi$ and~$g$ have the same order, the order of~$g$ divides~$\ind(A)$, and if we choose for~$g$ a generator of the cyclic group~$G_V$, we get that~$|G_V|$ divides~$\ind(A)$.

The more difficult part is to establish the converse. For this, let $\chi' \in \Ch(H)$ be an eigenvector for the left multiplication by~$\chi$ corresponding to the eigenvalue $\zeta \dim(V)$. Note that $\chi'$ is not necessarily the character of any module. By the Perron-Frobenius theorem, such an eigenvector is unique up to scalar multiples. Now, for any $\chi'' \in \Ch(H)$, 
$\chi' \chi''$ is also an eigenvector corresponding to this eigenvalue, and therefore there is a number~$\gamma(\chi'') \in K$\index{$\gamma$} such that
$\chi' \chi'' = \gamma(\chi'') \chi'$.
It is clear that
$$\gamma: \Ch(H) \rightarrow K$$
is an algebra homomorphism. Since $\Ch(H)$ is semisimple (cf.~\cite{Z2}, Lem.~2, p.~55), this shows that $K \chi'$ is a one-dimensional two-sided ideal of~$\Ch(H)$; in particular, $\chi'$ is central and we have
$$\zeta \dim(V) \chi' = \chi \chi' = \chi' \chi
= \gamma(\chi) \chi'$$
so that $\gamma(\chi) = \zeta \dim(V)$. Raising this equation to the $m$-th power, we get $\gamma(\chi^m) = \zeta^m \dim(V^{\o m})$. By decomposing
$V^{\o m}$ into simple modules, we get an equation of the form
$\chi^m = \sum_{i=1}^k k_i \chi_i$ for some nonnegative integers~$k_i$\index{$k_i$}.
Applying~$\gamma$, this equation becomes
$$\sum_{i=1}^k \zeta^m k_i n_i = \zeta^m \dim(V^{\o m}) = \gamma(\chi^m) = 
\sum_{i=1}^k k_i \gamma(\chi_i)$$
Since $\chi_i \chi' = \gamma(\chi_i) \chi'$, we have that 
$\gamma(\chi_i)$ is an eigenvalue of the left multiplication by~$\chi_i$.
But as discussed at the end of Paragraph~\ref{Index}, the absolute value of~$\gamma(\chi_i)$ is bounded by~$n_i=\dim(V_i)$. The above equality can therefore only hold if we have $\zeta^m k_i n_i = k_i \gamma(\chi_i)$ for all $i=1,\ldots,k$. This shows that we have $\gamma(\chi_i)=\zeta^m n_i$ if~$k_i \neq 0$, i.e., if $V_i$ appears as a constituent of~$V^{\o m}$. In particular, the numbers~$m$ for which $V_i$ appears as a constituent of~$V^{\o m}$ cannot be arbitrary, but can only appear in an $\ind(A)$-arithmetic progression\index{$\ind(A)$}.

Since~$J$, the intersection of the annihilators of all the tensor powers of~$V$, is zero, every simple module~$V_i$ appears as a constituent of some tensor power~$V^{\o m_i}$, and we therefore have 
$\gamma(\chi_i)=\zeta^{m_i} n_i$ for some number~$m_i$\index{$m_i$}, which is unique modulo~$\ind(A)$. Note that we can assume that $m_i=1$ whenever $V_i$ is a constituent of~$V$. Now define
$$g := \sum_{i=1}^k \zeta^{m_i} e_i$$
where $e_i \in Z(H)$\index{$e_i$} is the centrally primitive idempotent corresponding to~$V_i$. 
This is obviously a central element of order $\ind(A)$ that satisfies $\gamma(\chi_i)=\chi_i(g)$ for all $i=1,\ldots,k$. We claim that~$g$ is grouplike. For this, note that by construction~$g$ acts on~$V^{\o m}$ by multiplication with~$\zeta^m$. This implies that
both $\Delta(g)$ and~$g \o g$ act on~$V^{\o m} \o V^{\o l}$ by multiplication with~$\zeta^{m+l}$, so that $\Delta(g) - g \o g$
annihilates this module. Now the annihilator of the $H \o H$-module 
$V^{\o m} \o V^{\o l}$ is the sum $I_1 \o H + H \o I_2$, where 
$I_1$\index{$I_1$} is the annihilator of~$V^{\o m}$ and $I_2$\index{$I_2$} is the annihilator of~$V^{\o l}$. This shows that the intersection of all these annihilators
is zero, so that in particular $\Delta(g) - g \o g = 0$, i.e., $g$ is grouplike. Since~$g$ acts on~$V$ by multiplication with~$\zeta$, we have that~$g$ is an element of~$G_V$ of order~$\ind(A)$\index{$\ind(A)$}, so that conversely $\ind(A)$~divides~$|G_V|$. Furthermore, this shows that~$g$ generates~$G_V$.~\qed
\end{pf}
This proof also shows what the diagonal matrix~$D$ that appears in the second part of the Perron-Frobenius theorem is in the present case: In the situation and with the notation of the proof, we have 
$D=\diag(\zeta^{m_1},\ldots,\zeta^{m_k})$, so that~$D$ is the matrix representation of the map
$$\Ch(H) \rightarrow \Ch(H),~\varphi \mapsto (g \rightarrow \varphi)$$
with respect to the basis consisting of the irreducible characters. Here, the action appearing in this expression is defined by 
$(g \rightarrow \varphi)(h) = \varphi(hg)$. This holds since we have
$g \rightarrow \chi = \zeta \chi$ and therefore
$$g \rightarrow (\chi \varphi) = 
(g \rightarrow \chi) (g \rightarrow \varphi) 
= \zeta \chi (g \rightarrow \varphi)$$

If $V$ is simple, then every central grouplike element of~$H/J$, i.e., every grouplike element that is central in~$H/J$, acts on~$V$ by multiplication with a scalar. Therefore, $G_V = G(H/J) \cap Z(H/J)$ is exactly the set of central grouplike elements, and we get the following corollary:
\begin{corollary}
Suppose that~$V$ is simple. Then the group $G(H/J) \cap Z(H/J)$ of central grouplike elements of~$H/J$ is cyclic, and its order is equal to the index of imprimitivity of~$A$\index{$\ind(A)$}: 
$$\ind(A) = |G(H/J) \cap Z(H/J)|$$
\end{corollary}
As another consequence, we record the following relation between the three invariants that we have studied:
\begin{prop}
$\ind(A)$ divides $\exp(H)$ and $\ord(V)$.
\end{prop}
\begin{pf}
As before, let $J$ be the largest Hopf ideal contained in the annihilator of~$V$. That $\ind(A) = \exp(G_V)$ divides $\exp(H)$ follows from the fact that the exponent of Hopf subalgebras and quotients divide the exponent of the larger object (cf.~\cite{EG1}, Prop.~2.2.7, p.~132).
To see that $\ind(A)$ divides $\ord(V)$, note first that the order of~$V$ as a module over~$H$ is the same as its order as a module over~$H/J$, which implies that we can, as before, assume that~$J$ is zero. By the second part of the Perron-Frobenius theorem, we now have a diagonal matrix~$D$ such that $D A D^{-1} = \zeta A$, where~$\zeta$ is a primitive 
$\ind(A)$-th root of unity. If~$m$ is the order of~$V$, we know, for every irreducible character~$\chi_i$ of~$H$, that $\chi^m \chi_i$ contains the character $\chi_i$ at least $\mult(V)$ times, which shows that the diagonal entries of~$A^m$ are all strictly positive integers. By comparing the diagonal components on both sides of the equation~$D A^m D^{-1} = \zeta^m A^m$, we see that $\zeta^m=1$, which shows that $\ind(A)$ divides $m=\ord(V)$.
\qed
\end{pf}
One might conjecture that the index of~$A$ is actually always equal to the order of~$V$. However, this is not the case, as one can already see in the case of finite groups: Recall the example of a nonabelian group~$G$ of order~$pq$, which we have considered in Paragraph~\ref{Examplepq}, and the $p$-dimensional simple module~$V$ constructed there. First, the Hopf ideal~$J$ is zero in this case, because otherwise~$H/J$ would be the group ring of the factor group~$G(H/J)$, which would be necessarily cyclic, and so~$V$ would be a $p$-dimensional simple module of a cyclic group, which is impossible. Second, the center of~$G$ is trivial, as one can see by inspection or from the fact that a group that is cyclic modulo its center must be abelian. Therefore, the above corollary implies that $\ind(A)=1$\index{$\ind(A)$}. On the other hand, we have seen in Paragraph~\ref{Examplepq} that~$\ord(V)=3$ if~$p=5$ and~$q=11$.

\section{The Drinfel'd Double} \label{Sec:DrinfDouble}
\subsection[The Drinfel'd double]{} \label{DrinfDouble}
For an arbitrary finite-dimensional Hopf algebra~$H$, the Drinfel'd double is a Hopf algebra with underlying vector space~$H^* \o H$. The coalgebra structure is the tensor product coalgebra structure~$H^{* \cop} \o H$,
so that coproduct and counit are given by the formulas
$$\Delta(\varphi \o h) = (\varphi_\2 \o h_\1) \o (\varphi_\1 \o h_\2)
\qquad \varepsilon(\varphi \o h) = \varphi(1) \varepsilon(h)$$
The formula for the product is a little more involved; it reads
$$(\varphi \o h)(\varphi' \o h') = 
\varphi'_\1(S^{-1}(h_\3)) \varphi'_\3(h_\1) \;
\varphi \varphi'_\2 \o h_\2 h'$$
Finally, the antipode is given by the formula
$S(\varphi \o h) = (\varepsilon \o S(h))(S^{-1}(\varphi) \o 1)$.

One important point about the Drinfel'd double is that it is quasitriangular.
If $b_1,\ldots,b_n$ is a basis of~$H$ with dual basis $b^*_1,\ldots,b^*_n$,
then the R-matrix\index{$R$} is 
$$R = \sum_{i=1}^n (\varepsilon \o b_i) \o (b^*_i \o 1)$$
The associated Drinfel'd element~$u$\index{$u$} and its inverse are therefore 
$$u = \sum_{i=1}^n S^{-1}(b^*_i) \o b_i \qquad 
u^{-1} = \sum_{i=1}^n S^2(b^*_i) \o b_i$$
(cf.~\cite{M}, Chap.~10 for further details).

In the case of a semisimple Hopf algebra over a field of characteristic zero that we are considering, one can of course replace the inverse of the antipode by the antipode itself in the above formulas. In this case, the Drinfel'd double is also semisimple and cosemisimple (cf.~\cite{R2}, Prop.~7, p.~304), and if $\lambda \in H^*$\index{$\lambda$} is an integral that satisfies $\lambda(1)=1$, then
$\Gamma := \lambda \o \Lambda$ is an integral\index{$\Gamma$} of~$D(H)$ satisfying~$\varepsilon(\Gamma)=1$ (cf.~\cite{R2}, Thm.~4, p.~303). Furthermore, $u$ is a central element (cf.~\cite{M}, Prop.~10.1.4, p.~179) whose order is equal to the exponent of~$H$ (cf.~\cite{EG1}, Thm.~2.5, p.~133).

\subsection[Factorizability]{} \label{Factorize}
If $R_{21}$ arises from~$R$ by interchanging the tensorands, then the map
$$\Phi: D(H)^* \rightarrow D(H),~\psi \mapsto (\id \o \psi)(R_{21} R)$$
is bijective, i.e., the Drinfel'd double is factorizable 
(cf.~\cite{ResSem}, Def.~2.1, p.~543; see also \cite{SchneiderFact}, p.~1892). If restricted to the character ring of~$D(H)$, $\Phi$ induces an algebra isomorphism between the character ring and the center of the Drinfel'd double (cf.~\cite{DrinfAlmCocom}, Prop.~3.3, p.~327; see also \cite{SchneiderFact}, Thm.~2.1, p.~1892). It should be pointed out that,
under the identification 
$$H \o H^* \rightarrow D(H)^*,~h \o \varphi \mapsto 
(\varphi' \o h' \mapsto \varphi'(h) \varphi(h'))$$
the restriction of~$\Phi$ to the character ring is just the interchange of the tensorands. This holds since, for an element 
\mbox{$\eta=\sum_{k=1}^r h_k \o \varphi_k \in \Ch(D(H)) \subset H \o H^*$},
we have
\begin{align*}
\Phi(\eta) &= (\id \o \eta)(R_{21} R)
= \sum_{i,j=1}^n (b^*_j \o 1)(\varepsilon \o b_i) \;
\eta((\varepsilon \o b_j)(b^*_i \o 1)) \\
&= \sum_{i,j=1}^n (b^*_j \o b_i) \; \eta(b^*_i \o b_j)
=\sum_{k=1}^r \varphi_k \o h_k
\end{align*}
where we have used the fact that $\eta$ is a character in the third equality.

We are therefore in the following situation: Suppose that $\eta_1,\ldots,\eta_l$\index{$\eta_i$} are the different irreducible characters of~$D(H)$, with corresponding degrees $m_1,\ldots,m_l$\index{$m_i$}. We assume that $\eta_1=\varepsilon$ is the trivial character, so that $m_1=1$. The corresponding central characters, i.e., the algebra homomorphisms from the center~$Z(D(H))$ of the Drinfel'd double to~$K$, are then given by\index{$Z(D(H))$}\index{$\omega_i$}
$$\omega_i: Z(D(H)) \rightarrow K,~z \mapsto \frac{1}{m_i} \eta_i(z)$$
for $i=1,\ldots,l$\index{$\omega_i$}. This function takes the value~$1$ on the centrally primitive idempotent corresponding to~$\eta_i$ and the value~$0$ on the other centrally primitive idempotents. From the isomorphism between the character ring and the center, we get that the algebra homomorphisms from the character ring to the base field are the functions\index{$\xi_i$}
$$\xi_i: \Ch(D(H)) \rightarrow K,~\eta \mapsto \omega_i(\Phi(\eta))$$
for $i=1,\ldots,l$. 

As in Paragraph~\ref{Index}, we denote by~$\Ch_\Q(D(H))$\index{$\Ch_\Q(D(H))$} the $\Q$-algebra spanned by the irreducible characters~$\eta_1,\ldots,\eta_l$ over the rational numbers. We can say the following about the algebra structure:
\begin{prop}
For $d:=\exp(H)$, there are subfields $K_1,\ldots,K_s$ of the cyclotomic field~$\Q_d$ such that $\Ch_\Q(D(H)) \cong \bigoplus_{i=1}^s K_i$.
\end{prop}
\begin{pf}
The ordinary character ring~$\Ch(D(H))$ arises from the rational one by extension of scalars: $\Ch(D(H)) \cong K \o_\Q \Ch_\Q(D(H))$. From the fact that $\Ch(D(H))$ is isomorphic to~$Z(D(H))$, which is in turn isomorphic to 
$K^l$ as an algebra, we see that $\Ch(D(H))$ is commutative and semisimple.
(Alternatively, the quasitriangularity of~$D(H)$ implies that $\Ch(D(H))$
is commutative; the semisimplicity was already discussed in Paragraph~\ref{IndForm}.) As it cannot contain nilpotent elements, this implies that also 
$\Ch_\Q(D(H))$ is commutative and semisimple. By Wedderburn's theorem, we see that there are number fields $K_1,\ldots,K_s$ such that $\Ch_\Q(D(H)) \cong \bigoplus_{i=1}^s K_i$. We have to show that these fields can be embedded into~$\Q_d$. To see this, recall that 
$R_{21} R = (u \o u) \Delta(u^{-1})$ (cf.~\cite{M}, Thm.~10.1.13, p.~181) and therefore
$$m_i \xi_i(\eta_j) = m_i \omega_i(\Phi(\eta_j)) 
= (\eta_i \o \eta_j)(R_{21} R) 
= \omega_i(u) \omega_j(u) (\eta_i \eta_j)(u^{-1})$$
The product of the characters~$\eta_i$ and~$\eta_j$ can be written as a linear combination of the irreducible characters, i.e., we have
$\eta_i \eta_j = \sum_{k=1}^l N_{ij}^k \eta_k$ with coefficients~$N_{ij}^k$ that are integers. We therefore get
$$m_i \xi_i(\eta_j) 
= \omega_i(u) \omega_j(u) \sum_{k=1}^l N_{ij}^k m_k \omega_k(u)^{-1}$$
As we said above, $u$ has order~$d$, and therefore the elements~$\omega_k(u)$ are $d$-th roots of unity. This shows that
$\xi_i(\eta_j) \in \Q_d$, which implies that $\xi_i$ restricts to an algebra homomorphism from~$\Ch_\Q(D(H))$ to~$\Q_d$. Now, for a nonzero element $\eta \in K_j$, there is some~$\xi_i$ such that $\xi_i(\eta)$ is nonzero, and therefore the restriction of~$\xi_i$ is a nonzero homomorphism of $\Q$-algebras from~$K_j$ to~$\Q_d$, which must be an embedding of fields. 
\qed
\end{pf}

\subsection[The center of the character ring]{} \label{CentCharRing}
The Hopf algebra~$H$ can be considered as a subalgebra of the Drinfel'd double via the embedding $h \mapsto \varepsilon \o h$. If~$W$ is a module over the Drinfel'd double, we therefore get an $H$-module by restriction. However, this module has special properties: If $V$ is an arbitrary 
$H$-module, then the map
$$V \o W \rightarrow W \o V,~v \o w \mapsto 
\sum_{i=1}^n (b_i^* \o 1).w \o b_i.v$$
is an isomorphism of $H$-modules. This is a consequence of the condition
$$((1 \o h_\2) \o (1 \o h_\1)) R = R ((1 \o h_\1) \o (1 \o h_\2))$$
for the R-matrix, since we can write the above map as the composition of
the map $v \o w \mapsto R.(v \o w)$ with the usual isomorphism between
$V \o W$ and~$W \o V$. In fact, this observation can be extended to the statement that the modules over the Drinfel'd double form the so-called 
categorical center of the category of $H$-modules (cf.~\cite{Kas}, Thm.~XIII.5.1, p.~333).

We therefore have an algebra homomorphism\index{$\Res$}
$$\Res: \Ch(D(H)) \rightarrow Z(\Ch(H))$$
from the character ring of the Drinfel'd double to the center of the character ring of~$H$, which assigns to every character its restriction to~$H$. The basic property of this map is the following:
\begin{prop}
$\Res$ is surjective.
\end{prop}
\begin{pf}
From the trivial $H$-module~$K$, we get the induced module 
$D(H) \o_H K$ over~$D(H)$. The vector space isomorphism 
$H^* \rightarrow D(H) \o_H K,~\varphi \mapsto (\varphi \o 1) \o_H 1$
becomes an isomorphism of $D(H)$-modules if we define the $D(H)$-action on~$H^*$ by the formula 
$(\varphi \o h).\varphi' := \varphi(h.\varphi')$,
where 
$$h.\varphi = \varphi_\1(S(h_\2)) \varphi_\3(h_\1) \varphi_\2$$
Now suppose that $\chi \in \Ch(H)$. Then we have
\begin{align*}
h.(\varphi \chi) &= 
\varphi_\1(S(h_\4)) \chi_\1(S(h_\3)) \varphi_\3(h_\1) \chi_\3(h_\2) 
\; \varphi_\2 \chi_\2 \\
&= \varphi_\1(S(h_\4)) \chi_\2(S(h_\3)) \varphi_\3(h_\1) \chi_\1(h_\2) 
\; \varphi_\2 \chi_\3 \\
&= \varphi_\1(S(h_\2)) \varphi_\3(h_\1) \; \varphi_\2 \chi
= (h.\varphi) \chi
\end{align*}
From this, we see that the map
$$\Ch(H) \rightarrow \End_{D(H)}(H^*),~\chi \mapsto 
(\varphi \mapsto \varphi \chi)$$
is an algebra antihomomorphism. It is even an anti-isomorphism, since every $D(H)$-linear endomorphism of~$H^*$ must be in particular $H^*$-linear, and therefore be given by right multiplication with some $\chi \in H^*$. But the condition
$h.\chi = (h.\varepsilon) \chi = \varepsilon(h) \chi$
requires that $\chi \in \Ch(H)$.

Since $\Ch(H)$ is isomorphic to the endomorphism ring of a semisimple module, it is itself semisimple (cf.~\cite{FarbDennis}, Prop.~1.8, p.~36). In particular, if $p \in \Ch(H)$ is a centrally primitive idempotent, then the right multiplication by~$p$ on~$H^*$ is the projection onto some isotypical component of the $D(H)$-module~$H^*$, i.e., in a decomposition of~$H^*$ into a direct sum of simple $D(H)$-modules, it is the projection onto the sum of all simple modules of a given isomorphism type. If $e \in Z(D(H))$ is the centrally primitive idempotent that corresponds to this isomorphism type, then $e$ also acts as the projection onto this isotypical component, so that $e.\varphi = \varphi p$ for all $\varphi \in H^*$.

Now let $\eta:=\Phi^{-1}(e) \in \Ch(D(H))$. If we identify $D(H)^*$ 
with~$H\o H^*$ as described in Paragraph~\ref{Factorize}, we can write 
$\eta = \sum_{k=1}^r h_k \o \varphi_k$. As discussed there, we then have
$e = \sum_{k=1}^r \varphi_k \o h_k$ and therefore
$$p = \varepsilon p = e.\varepsilon 
= \sum_{k=1}^r \varphi_k \varepsilon(h_k)$$
which gives 
$p(h) = \sum_{k=1}^r \varepsilon(h_k) \varphi_k(h) 
= \eta(\varepsilon \o h)$. This shows that $\Res(\eta)=p$, which establishes the claim since the centrally primitive idempotents form a basis of~$Z(\Ch(H))$.
\qed
\end{pf}
We note that the isomorphism between~$\Ch(H)$ and~$\End_{D(H)}(H^*)$, which is crucial in the above proof, is dual to the isomorphism constructed in~\cite{Z3}, Thm.~1, p.~2849. It gives another proof of the result that the character ring of a semisimple cosemisimple Hopf algebra is itself semisimple (cf.~\cite{Z2}, Lem.~2, p.~55; see also~\cite{CohenZhu}, Thm.~2.4, p.~164), which also has the advantage that it does not really rely on our assumption that the base field has characteristic zero (cf.~\cite{CohenGenChar}, Thm.~3.6, p.~62).

From the preceding proposition, we can learn something about the arithmetic of the character ring of~$H$. In contrast to the case of the Drinfel'd double, the span~$\Ch_\Q(H)$\index{$\Ch_\Q(H)$} of the irreducible characters of~$H$ over~$\Q$ is, in general, not commutative, but it is, as we just saw, still semisimple. The Wedderburn theorem therefore tells that it is isomorphic to a direct sums of matrix rings over division algebras. The center of~$\Ch_\Q(H)$ is then the sum of the centers of these division algebras, which are some number fields about which we can now say the following:
\begin{thm}
For $d:=\exp(H)$, there are subfields $K_1,\ldots,K_t$ of the cyclotomic field~$\Q_d$ such that $Z(\Ch_\Q(H)) \cong \bigoplus_{i=1}^t K_i$.
\end{thm}
\begin{pf}
Since forming centers is compatible with the extension of scalars, we get from $\Ch(H) \cong K \o_\Q \Ch_\Q(H)$ that also
$Z(\Ch(H)) \cong K \o_\Q Z(\Ch_\Q(H))$ (cf.~\cite{FarbDennis}, Chap.~3, Exerc.~22, p.~103). Therefore, the preceding proposition implies that $\Res$ restricts to a surjective algebra homomorphism from $\Ch_\Q(D(H))$ to~$Z(\Ch_\Q(H))$. Since a surjective algebra homomorphism between semisimple algebras is always the projection onto some two-sided ideals, the assertion follows from Proposition~\ref{Factorize}.
\qed
\end{pf}
In the case where $H=K^G$ is the ring of functions on a finite group~$G$, we have explained in Paragraph~\ref{OrderMult} that the character ring of~$H$ is the group ring~$K[G]$. In this case, the above result tells that the center of the group ring~$\Q[G]$ over the rational numbers is isomorphic to a product of subfields of the cyclotomic fields determined by the exponent. In the case of groups, there is a stronger result, which is due to R.~Brauer: This cyclotomic field is actually a splitting field for the group ring (cf.~\cite{Isaacs}, Thm.~(10.3), p.~161; \cite{Serre1}, Sec.~12.3, Thm.~24, p.~94). One may conjecture that the analogous result holds for semisimple Hopf algebras: If $d=\exp(H)$, the field~$\Q_d$ is a splitting field for~$\Ch_\Q(H)$.

\subsection[The third formula for the Frobenius-Schur indicators]{} \label{ThirdForm}
The induced module $D(H) \o_H K$ of the trivial module played an important role in the proof of Proposition~\ref{CentCharRing}. We now work out the character of an arbitrary induced module, which will lead to a third formula for the Frobenius-Schur indicators. From the adjoint action of the integral~$\Gamma$ of~$D(H)$, we get the projection\index{$P$}
$$P: D(H) \rightarrow Z(D(H)),~x \mapsto \Gamma_\1 x S(\Gamma_\2)$$
to the center of the Drinfel'd double. It is obvious that this map is the identity on the center; the fact that it really maps to the center follows from~\cite{LR2}, Lem.~1.2, p.~270.

Now suppose that $V$ is an arbitrary $H$-module with character~$\chi$. Consider the induced module $D(H) \o_H V$, and denote its character by~$\eta$. Using the identification between $H \o H^*$ and $D(H)^*$ described in Paragraph~\ref{Factorize}, we have the following formula for this character:
\begin{prop}
For all $x \in D(H)$, we have 
$\eta(x) = \dim(H)\, (\Lambda \o \chi)(P(x))$.
\end{prop}
\begin{pf}
The map 
$H^* \o V \rightarrow D(H) \o_H V,~\varphi \o v \mapsto 
(\varphi \o 1) \o_H v$
is a $K$-linear isomorphism. We can therefore endow~$H^* \o V$ with a uniquely determined $D(H)$-module structure such that this map is $D(H)$-linear. Explicitly, this module structure is given as
$$(\varphi \o h).(\varphi' \o v) = 
\varphi'_\1(S(h_\3)) \varphi'_\3(h_\1) \; \varphi \varphi'_\2 \o h_\2.v$$
Suppose first that $x$ is central, and write it in the form
$x = \sum_{i=1}^r \varphi_i \o h_i$. We then have
$$x.(\varphi \o v) = x (\varphi \o 1).(\varepsilon \o v)
= (\varphi \o 1) x.(\varepsilon \o v) 
= \sum_{i=1}^r \varphi \varphi_i \o h_i.v$$
This shows that $x$ acts as the sum of the tensor products of the right multiplications by~$\varphi_i$ and the module actions of~$h_i$. Since the trace of the right multiplication by~$\varphi_i$ is 
$\dim(H) \, \varphi_i(\Lambda)$ (cf.~\cite{LR2}, Prop.~2.4, p.~273) and the trace of the module action of~$h_i$ is
$\chi(h_i)$, we get that
$$\eta(x) = \dim(H)\sum_{i=1}^r \varphi_i(\Lambda) \chi(h_i)
= \dim(H) \, (\Lambda \o \chi)(x)$$
which proves the assertion in this case since $P(x)=x$.

If $x$ is not necessarily central, we argue as follows: Since $\eta$ is a character, we have 
$$\eta(P(x)) = \eta(\Gamma_\1 x S(\Gamma_\2)) 
= \eta(x S(\Gamma_\2) \Gamma_\1) = \eta(x)$$ 
Since $P(x)$ is central, this yields
$\eta(x) = \eta(P(x)) = \dim(H) \, (\Lambda \o \chi)(P(x))$.
\qed
\end{pf}
As a consequence, we get a third formula for the Frobenius-Schur indicators:
\begin{corollary}
$\displaystyle \nu_m(\chi)=\frac{1}{\dim(H)} \eta(u^m)$
\end{corollary}
\begin{pf}
As discussed in Paragraph~\ref{DrinfDouble}, the Drinfel'd element
is given by the formula
$u = \sum_{i=1}^n S(b^*_i) \o b_i$,
where $b_1,\ldots,b_n$ is a basis of~$H$ with dual basis $b^*_1,\ldots,b^*_n$. Since $u$ is central, we then have
\begin{align*}
u^2 &= \sum_{i=1}^n (S(b^*_i) \o 1) u (\varepsilon \o b_i) \\
&= \sum_{i_1,i_2=1}^n S(b^*_{i_2}) S(b^*_{i_1}) \o b_{i_1} b_{i_2}
= \sum_{i_1,i_2=1}^n S(b^*_{i_1} b^*_{i_2}) \o b_{i_1} b_{i_2}
\end{align*}
Repeating this argument, we get
$$u^m = \sum_{i_1,\ldots,i_m=1}^n S(b^*_{i_1} b^*_{i_2} \ldots  
b^*_{i_m}) \o b_{i_1} b_{i_2} \ldots b_{i_m}$$
Since $P(u^m)=u^m$, we get from the preceding proposition that
\begin{align*}
\eta(u^m) &= \dim(H) \sum_{i_1,\ldots,i_m=1}^n 
S(b^*_{i_1} b^*_{i_2} \ldots b^*_{i_m})(\Lambda) \;
\chi(b_{i_1} b_{i_2} \ldots  b_{i_m})\\
&= \dim(H) \, \chi(\Lambda_\1 \Lambda_\2 \ldots \Lambda_\m)
= \dim(H) \, \nu_m(\chi)
\end{align*}
where we have used that $S(\Lambda) = \Lambda$. 
\qed
\end{pf}

\section{Examples} \label{Sec:Examples}
\subsection[The setup]{} \label{SetupCocentAbExt}
In this section, we illustrate the theory developed so far by considering the example of a very special class of Hopf algebras, which are, however, in general neither group rings nor dual group rings. We continue to work over an algebraically closed field~$K$ of characteristic zero. Suppose that $F$\index{$F$} and $G$ are two finite groups and
that $F$ acts on~$G$ by group automorphisms. We denote this action
by
$$F \times G \rightarrow G,~(x,g) \mapsto x.g$$
Consider the ring~$K^G$ of functions on~$G$, and denote by~$b_g$\index{$b_g$} the function that takes the value~$1$ on the element $g \in G$ and the value~$0$ on all other group elements; these functions form a basis of~$K^G$. The action of~$F$ on~$G$ leads to an action of~$F$ on~$K^G$ by defining
$$x.b_g = b_{x.g}$$

On the vector space $K^G \o K[F]$, we introduce a Hopf algebra
structure as follows: For the coalgebra structure, we take the usual tensor product of the coalgebras~$K^G$ and~$K[F]$, so that the coproduct and the counit are given on basis elements by the formulas
$$\Delta(b_g \o x) = \sum_{\substack{g_1, g_2 \in G\\g_1g_2=g}}
(b_{g_1} \o x) \o (b_{g_2} \o x) \qquad
\varepsilon(b_g \o x) = \delta_{g,1}$$
for $g \in G$ and $x \in F$. For the algebra structure, we take the smash product with respect to the action described above, so that the product is given on basis elements by the formula
$$(b_{g_1} \o x_1)(b_{g_2} \o x_2) 
= \delta_{g_1,x_1.g_2} \; b_{g_1} \o x_1 x_2$$
for $g_1, g_2 \in G$ and $x_1, x_2 \in F$, and the unit element is the 
tensor product of the unit elements, which takes the form
$1 = \sum_{g \in G} b_g \o 1$
if expanded in terms of basis elements. Finally, we introduce the antipode by requiring that we have
$$S(b_g \o x) = (1 \o x^{-1})(b_{g^{-1}} \o 1)$$
on basis elements.

Besides verifying directly that these structure elements turn
$K^G \o K[F]$ into a Hopf algebra, it is possible to view them as arising
from at least two constructions. First, one can see $K^G$ as a Yetter-Drinfel'd Hopf algebra over~$K[F]$: The action of~$K[F]$ on~$K^G$ is the one introduced above and the coaction of~$K[F]$ on~$K^G$ is trivial. The Hopf algebra structure on~$K^G \o K[F]$ introduced above is then the associated Radford biproduct (cf.~\cite{M}, Thm.~10.6.5, p.~209).

The second, slightly more standard view of this Hopf algebra is to consider it as an extension of~$K[F]$ by~$K^G$: With respect to the above action of~$K[F]$ on~$K^G$ and the trivial coaction of~$K^G$ on~$K[F]$, $K[F]$ and~$K^G$ form an abelian matched pair in the sense of~\cite{Hofst2}, Def.~1.1, p.~266, and the above Hopf algebra structure is the one associated there with the trivial cocycles (cf.~\cite{Hofst2}, Sec.~3, p.~273). This ties the considerations below with the investigations in~\cite{KMM}, where the second indicator was analyzed in the more general case where nontrivial cocycles are present.

It should be noted that the Drinfel'd double of~$K[F]$, which we have considered in Paragraph~\ref{DrinfDouble}, reappears as a special case of this construction. For the Drinfel'd double, we have $G=F^{\op}$, the opposite group of~$F$ in which the multiplication is reversed, and the action is just conjugation, i.e., we have $x.g := xgx^{-1}$, where the product is taken in~$F$, and not in~$G=F^{\op}$.

\subsection[The coefficients]{} \label{Coeff}
The element $\Lambda := {{\scriptstyle 1/}{}_{|F|}}\sum_{x \in F} b_1 \o x$ is an integral of~$K^G \o K[F]$ such that $\varepsilon(\Lambda) = 1$ (cf.~\cite{KMM}, Sec.~4, p.~896). In order to describe
its Sweedler powers, we introduce some notation. Suppose that~$m$
and~$k$ are natural numbers that are relatively prime. For an arbitrary
integer~$j$, we denote by~$\l j\r \in
\{0,1,\ldots,m-1\}$\index{$\l j\r$} the remainder
under division by~$m$. For $g \in G$ and $y \in F$, we define the set\index{$G_{m,k}(g,y)$}
$$G_{m,k}(g,y) := \{x \in F \mid x^m = y \text{ and }
\prod_{j=0}^{m-1} x^{-\l j/k\r}.g  = 1 \}$$ 
and denote the cardinality of this set by $z_{m,k}(g,y)$\index{$z_{m,k}(g,y)$}. Here, the quotient $j/k$ should be understood modulo~$m$, so that, if $l$ is an integer that satisfies 
$kl \equiv 1 \pmod{m}$, we have $\l j/k\r = \l jl\r$. We will see in Paragraph~\ref{SweedlPowIntExa} how these numbers relate to the Sweedler powers of~$\Lambda$; before we discuss this, however, let us put down their basic properties:
\begin{prop}
If $m$, $k$, and $q$ are pairwise relatively prime natural numbers, the following equations hold for $x, y \in F$ and $g \in G$:
\begin{enumerate}
\item $z_{m,k}(x.g,y) = z_{m,k}(g,x^{-1}yx)$
\item $z_{m,k}(g,y) = 0$ unless $y.g = g$
\item $z_{mk,q}(g,y^k) = z_{m,q}(g,y)$ if $k$ is relatively prime to~$|F|$ and~$|G|$
\item $z_{m,kq}(g,y^k) = z_{m,q}(g,y)$ if $k$ is relatively prime to~$|F|$
\end{enumerate} 
\end{prop}
\begin{pf}
To prove the first statement, note that $z \in G_{m,k}(g,x^{-1}yx)$ if and only if
$$z^m=x^{-1}yx \quad \text{and} \quad 
\prod_{j=0}^{m-1} z^{-\l jl\r}.g = 1$$
where $l$ is an integer that satisfies $kl \equiv 1 \pmod{m}$. Acting with~$x$ on the second equation yields
$$1 = \prod_{j=0}^{m-1} x.(z^{-\l jl\r}.g)= 
\prod_{j=0}^{m-1} (xzx^{-1})^{-\l jl\r}x.g$$
so that $z \in G_{m,k}(g,x^{-1}yx)$ if and only if $xzx^{-1} \in G_{m,k}(x.g,y)$. This proves the first statement.

Since the element~$x$ does not appear in the remaining statements, we will in the rest of the proof use it again as a free variable. It is clear from the definition that the number~$z_{m,k}(g,y)$ depends only on the residue of~$k$ modulo~$m$, and we therefore can assume that~$k \in I_{m-1}$. To prove the second statement, we first note the following: If $l$ is again an integer that satisfies 
$kl \equiv 1 \pmod{m}$, we have for $j=1,\ldots,k-1$ that
$$\l (j+m-k)l\r +1 = \l jl\r$$
For we certainly have $(j+m-k)l+1 \equiv jl \pmod{m}$, and in addition we have
$\l (j+m-k)l\r +1 \le m-1$, since the value~$m$ is attained for~$j=0$.

Likewise, we have for $j=k,\ldots,m-1$ that
$$\l (j-k)l\r +1 = \l jl\r$$
since we obviously have $(j-k)l+1 \equiv jl \pmod{m}$, but furthermore also have $\l (j-k)l\r +1 \le  m-1$, because
$\l (j-k)l\r +1 = m$ holds for $j=m$.

Suppose now that $x \in G_{m,k}(g,y)$. By definition, we then have 
$x^m = y$ and $\prod_{j=0}^{m-1} x^{-\l jl\r}. g  = 1$. Applying~$x^{-1}$ yields
$$1 = \prod_{j=0}^{m-1} x^{-\l jl\r-1}.g
= (\prod_{j=0}^{m-k} x^{-\l jl\r-1}.g) 
(\prod_{j=m-k+1}^{m-1} x^{-\l jl\r-1}.g)$$
The last equation shows that the two factors are inverses of each other, and we can therefore interchange the factors to obtain
\begin{align*}
1 &= (\prod_{j=m-k+1}^{m-1} x^{-\l jl\r-1}.g)
(\prod_{j=0}^{m-k} x^{-\l jl\r-1}.g) \\
&= (\prod_{j=1}^{k-1} x^{-\l (j+m-k)l\r-1}.g)
(\prod_{j=k}^{m} x^{-\l (j-k)l\r-1}.g)
\end{align*}
By the considerations above, we can rewrite this in the form
$$1 = (\prod_{j=1}^{k-1} x^{-\l jl\r}.g)(\prod_{j=k}^{m-1} x^{-\l jl\r}.g)
x^{-m}.g = (\prod_{j=1}^{m-1} x^{-\l jl\r}.g) x^{-m}.g$$
Comparing this with the assumption, we arrive at $x^{-m}.g=g$, which proves the second statement. Note that this shows that~$x \in G_{m,k}(g,y)$ if and only if 
$$x^m=y \qquad x^m.g=g \qquad \prod_{j=0}^{m-1} x^{-jl}.g = 1$$

For the third statement, we claim that the stronger statement
$G_{mk,q}(g,y^k) = G_{m,q}(g,y)$ holds. Since~$mk$ and~$q$ are relatively prime, we can find an integer~$r$ such that 
$qr \equiv 1 \pmod{mk}$. By the characterization of these sets just obtained, we have $x \in G_{mk,q}(g,y^k)$ if and only if
$$x^{mk}=y^k \qquad x^{mk}.g=g \qquad \prod_{j=0}^{mk-1} x^{-jr}.g=1$$
Since $k$ is relatively prime to~$|F|$, the first and the second condition can be written in the form $x^m=y$, resp.~$x^m.g=g$. Therefore, the third condition can be written in the form
$$1 = \prod_{i=0}^{k-1} \prod_{j=0}^{m-1} x^{-(mi+j)r}.g 
= \prod_{i=0}^{k-1} \prod_{j=0}^{m-1} x^{-jr}.g 
= (\prod_{j=0}^{m-1} x^{-jr}.g)^k$$
Since $k$ is also relatively prime to~$|G|$, this shows that 
$x \in G_{mk,q}(g,y^k)$ if and only if
$$x^m=y \qquad x^m.g=g \qquad \prod_{j=0}^{m-1} x^{-jr}.g=1$$
which, since also $qr \equiv 1 \pmod{m}$, means that $x \in G_{m,q}(g,y)$.

For the fourth statement, observe that the assumptions assure that $m$ and~$kq$ are relatively prime, so that the number $z_{m,kq}(g,y^k)$ is actually defined and we can find an integer~$l$ such that
$kql \equiv 1 \pmod{m}$. We have that $x \in G_{m,q}(g,y)$ if and only if 
$x^k \in G_{m,kq}(g,y^k)$, since $x \in G_{m,q}(g,y)$ if and only if
$$x^m=y \qquad x^m.g=g \qquad \prod_{j=0}^{m-1} x^{-jkl}.g=1$$
and $x^k \in G_{m,kq}(g,y^k)$ if and only if
$$x^{mk}=y^k \qquad x^{mk}.g=g \qquad \prod_{j=0}^{m-1} (x^k)^{-jl}.g=1$$
This shows that $z_{m,kq}(g,y^k) = z_{m,q}(g,y)$, as asserted.
\qed
\end{pf}
Without the assumption that $k$ is relatively prime to~$|F|$, the proof of the fourth statement above still shows that $x^k \in G_{m,kq}(g,y^k)$ if 
$x \in G_{m,q}(g,y)$. By setting $q=1$ and~$y=1$, we therefore see that the map
$$G_{m,1}(g,1) \rightarrow G_{m,k}(g,1),~x \mapsto x^k$$
is well-defined. If, as above, $l$ is an integer that satisfies~$kl \equiv 1 \pmod{m}$, we see by setting~$y=1$ and 
replacing~$k$ by~$l$ and~$q$ by~$k$ that the map
$$G_{m,k}(g,1) \rightarrow G_{m,1}(g,1),~x \mapsto x^l$$
is well-defined. Since $x^m=1$ implies that $x^{kl}=x$, these maps are inverses of each other, so that we have the following corollary:
\begin{corollary}
If $m$ and $k$ are relatively prime natural numbers, we have for $y \in F$ and 
$g \in G$ that
$$z_{m,k}(1,y) = z_{m,1}(1,y)  \qquad  z_{m,k}(g,1) = z_{m,1}(g,1)$$
\end{corollary}
Here, the first assertion is obvious, since the condition 
$\prod_{j=0}^{m-1} x^{-\l j/k\r}.g=1$ is always satisfied if $g=1$.

\subsection[Sweedler powers of the integral]{} \label{SweedlPowIntExa}
The relation of the numbers~$z_{m,k}(g,y)$ with the Sweedler powers of~$\Lambda$ now is that they are essentially the coefficients in the expansion of the Sweedler powers in terms of the basis elements:
\begin{prop}
For relatively prime natural numbers~$m$ and~$k$, we have
$$\Lambda^{[m,k]} 
= \frac{1}{|F|} \sum_{g \in G,\, y \in F} z_{m,k}(g,y) \; b_g \o y$$
\end{prop}
\begin{pf}
Since
$$\Lambda_\1 \o \ldots \o \Lambda_\m =
\frac{1}{|F|} \sum_{x \in F} \sum_{\substack{g_1,g_2, \ldots ,g_m \in G \\ g_1 g_2 \ldots g_{m-1} g_m =1}} 
(b_{g_1} \o x) \o (b_{g_2} \o x) \o \ldots \o (b_{g_m} \o x)$$ 
we get  for the Sweedler power~$\Lambda^{[m,k]}$ the 
expression
\begin{align*}
&\frac{1}{|F|} \sum_{x \in F} \sum_{\substack{g_1,g_2, \ldots ,g_m \in G \\
g_1 g_2 \ldots g_{m-1} g_m  = 1}} (b_{g_{1}} \o x)
(b_{g_{1+\l k\r}} \o x) (b_{g_{1+\l 2k\r}} \o x)
\ldots (b_{g_{1+\l (m-1)k\r}} \o x) \\
&=\frac{1}{|F|} \sum_{x \in F} \sum_{\substack{g_1,g_2, \ldots ,g_m \in G \\
g_1 g_2 \ldots g_{m-1} g_m = 1}} 
(b_{g_1} b_{x.g_{1+\l k\r}} b_{x^2.g_{1+\l 2k\r}} \ldots b_{x^{m-1}.g_{1+\l (m-1)k\r}})
\o x^m
\end{align*}
If $l$ is an integer that satisfies $kl \equiv 1 \pmod{m}$, we can rearrange the product in the first tensorand above as follows:
$$\prod_{j=0}^{m-1} b_{x^j.g_{1+\l jk\r}}
= \prod_{j=0}^{m-1} b_{x^{\l lj\r}.g_{1+\l jlk\r}} =
\prod_{j=0}^{m-1} b_{x^{\l lj\r}.g_{1+j}}$$ 
Since a product of such basis elements is nonzero only if all the factors are equal, this expression is nonzero only if $g_{1+j} = x^{-\l lj\r}.g$ for
$j=0,\ldots,m-1$, where $g:=g_1$. With this notation, the condition 
$g_1 g_2 \ldots g_{m-1} g_m =1$ takes the form
$\prod_{j=0}^{m-1} x^{-\l lj\r}.g = 1$. We can therefore write the
Sweedler power in the form
\begin{align*}
\Lambda^{[m,k]}
=\frac{1}{|F|} \sum_{\substack{x \in F,\, g \in G\\
\prod_{j=0}^{m-1} x^{-\l lj\r}.g = 1}} b_g \o x^m
\end{align*}
The stated formula now follows by collecting the terms.
\qed
\end{pf}
In combination with Proposition~\ref{SweedlPowCent}, this formula can be used to give an alternative proof of the first two assertions in Proposition~\ref{Coeff}: Since $\Lambda^{[m,k]}$ is central, we have
$$(1 \o x)\Lambda^{[m,k]} = \Lambda^{[m,k]}(1 \o x) \qquad 
(b_g \o 1)\Lambda^{[m,k]} = \Lambda^{[m,k]}(b_g \o 1)$$
for $x \in F$ and $g \in G$, and by writing down these two equations in terms of the coefficients~$z_{m,k}(g,y)$, we get the first two assertions in Proposition~\ref{Coeff}, as the reader is invited to check. A consequence of the third and the fourth assertion can be deduced from Corollary~\ref{SweedlPowInt}: Under the stronger additional assumption that both~$k$ and~$q$ are relatively prime to $\dim(H)=|F||G|$, we have
$$\Lambda^{[mk,q]} = \Lambda^{[mkq]} = \Lambda^{[m,kq]}$$
and therefore $z_{mk,q}(g,y) = z_{m,kq}(g,y)$ for all $y \in F$ and 
all~$g \in G$.

\subsection[The irreducible modules]{} \label{IrredMod}
Our next goal is to find an explicit formula for the higher Frobenius-Schur indicators of these Hopf algebras. For this, we need a description of the irreducible modules and their characters, which we take from~\cite{KMM}, Sec.~3. Fix an element~$g \in G$ and denote the
stabilizer of~$g$ by~$F_g$\index{$F_g$}. The simple 
$K^G \o K[F]$-modules arise from simple $K[F_g]$-modules as follows:
\begin{prop}
For a $K[F_g]$-module~$W$, the induced module
$V:=K[F] \o_{K[F_g]} W$
becomes a $K^G \o K[F]$-module via
$$(b_h \o x).(y \o_{K[F_g]} w):= \delta_{h,xy.g} \; xy \o_{K[F_g]} w$$
If $W$ is a simple $K[F_g]$-module, then~$V$ is a simple $K^G \o K[F]$-module. Moreover, every simple $K^G \o K[F]$-module is isomorphic to a module that arises in this way.
\end{prop}
A proof of this proposition can be found in~\cite{KMM}, Cor.~3.5, p.~895.
To understand this module structure, it may be helpful to consider its restriction to the subalgebras that correspond to the two tensor factors, where the defining formula takes the simpler forms
$$(1 \o x).(y \o_{K[F_g]} w) =  xy \o_{K[F_g]} w \qquad
(b_h \o 1).(y \o_{K[F_g]} w) = \delta_{h,y.g} \; y \o_{K[F_g]} w$$
Note also that $b_h \o x$ acts identically as zero if~$g$ and~$h$ are not in the same orbit under the action of $F$. Furthermore, if
$z_1,\ldots,z_n$\index{$z_i$} is a system of representatives of the cosets $z F_g \in F/F_g$, every element~$v$ in~$V$ can be expressed uniquely in the form
$$v = \sum_{i=1}^n z_i \o_{K[F_g]} w_i$$
The action of~$1 \o x$ preserves the space~$z_i \o_{K[F_g]} W$ if and
only if $xz_i \in z_i F_g$, i.e., if $z_i^{-1} xz_i \in F_g$,
which means that~$x \in F_{z_i.g}$. If we use this decomposition to compute the character~$\chi$ of~$V$, we get the formula
\begin{align*}
\chi(b_h \o x) &= \chi((1 \o x)(b_h \o 1)) \\
&= \sum_{\substack{i=1 \\ x \in F_{z_i.g}}}^n \delta_{h,z_i.g}\eta(z_i^{-1}xz_i)
= \sum_{\substack{i=1 \\ x \in F_h}}^n \delta_{h,z_i.g}\eta(z_i^{-1}xz_i)\\
&=
\begin{cases}
\displaystyle
\frac{1}{|F_g|}\sum_{\substack{z \in F\\h=z.g}}\eta(z^{-1}xz)&: x \in F_h\\
0&: x \notin F_h
\end{cases}
\end{align*}
where~$\eta$ is the character of~$W$. Combining this formula with Proposition~\ref{SweedlPowIntExa}, we can now calculate the indicators of this character:
\begin{corollary}
For relatively prime natural numbers~$m$ and~$k$, we have
$$\chi(\Lambda^{[m,k]}) =
\frac{1}{|F_g|} \sum_{y \in F_{g}} z_{m,k}(g,y) \; \eta(y)$$
\end{corollary}
\begin{pf}
This follows by inserting the preceding results:
\begin{align*}
&\chi(\Lambda^{[m,k]})
= \frac{1}{|F|} \sum_{h \in G,\, x \in F} z_{m,k}(h,x) \; \chi(b_h \o x)\\
&= \frac{1}{|F|} \sum_{h \in G,\, x \in F_h} z_{m,k}(h,x) \; \chi(b_h \o x)
= \frac{1}{|F_g||F|} \sum_{z \in F,\, x \in F_{z.g}} z_{m,k}(z.g,x) \; \eta(z^{-1}xz) \\
&= \frac{1}{|F_g||F|} \sum_{z \in F,\, x \in F_{z.g}} z_{m,k}(g,z^{-1}xz) \; \eta(z^{-1}xz) 
= \frac{1}{|F_g|} \sum_{y \in F_{g}} z_{m,k}(g,y) \; \eta(y)
\end{align*}
Here, we have used Proposition~\ref{Coeff} for the fourth equality.
\qed
\end{pf}

\subsection[Nonintegral indicators]{} \label{ExampleA4}
We have now accumulated enough facts about this class of examples to employ it for the refutation of several conjectures that are suggested by previous considerations. To begin, recall that we have seen in Corollary~\ref{SecondForm} that the Frobenius-Schur indicators are integers if the exponent is squarefree. Let us give an example that they are not integers in general: Let $G=A_4$\index{$A_4$}, the alternating group on four letters, and let $F$ be cyclic of order~$9$. Choose a generator~$c$ of~$F$\index{$c$}. An action of~$F$ on~$G$ is then completely
described by the action of the generator~$c$; we require that this
action be the conjugation by the 3-cycle $\tau:=(1,2,3)$, so that
$$c.\sigma = \tau \sigma \tau^{-1}$$
We want to discuss the third Sweedler powers of the integral, so that $m=3$.
Consider the element $g=(1,4,3)$, and look at the equations
\mbox{$\prod_{j=0}^{m-1} c^{-\l j/k\r}.g = 1$}, where, as explained in Paragraph~\ref{Coeff}, the quotient~$j/k$ should be understood modulo~$m$.
In the case $k=1$ we have
$$
g(c^{-1}.g)(c^{-2}.g) = (1,4,3)(3,4,2)(2,4,1) = \id
$$
and in the case $k=2$ we have
$$
g(c^{-2}.g)(c^{-1}.g) = (1,4,3)(2,4,1)(3,4,2) \neq \id
$$
Using this, we can determine the sets $G_{3,k}(g,c^3)$. Note first that these sets are subsets of
$\{x \in F \mid x^3=c^3\}=\{c,c^4,c^7\}$, and, since $c^3.g=g$, they will
either coincide with this set or be empty. Now the above equations show that $G_{3,1}(g,c^3)$ coincides with
this set, whereas $G_{3,2}(g,c^3)=\emptyset$. This shows that
$$z_{3,1}(g,c^3) = 3 \qquad z_{3,2}(g,c^3)= 0$$
By the fourth assertion in Proposition~\ref{Coeff}, this implies
$$z_{3,1}(g,c^6) = 0 \qquad z_{3,2}(g,c^6)= 3$$
which can, of course, also be verified directly.
Furthermore, we have $z_{m,1}(g,1)=z_{m,2}(g,1)=3$.

Now let $\zeta$ be a primitive third root of unity, and let $\eta$ be the irreducible character
of $F_{g} =\left\langle c^3\right\rangle$ defined by $\eta (c^{3i})=\zeta^i$. As described in Paragraph~\ref{IrredMod},
we can associate with the corresponding one-dimensional $F_g$-module a three-dimensional $K^G \o K[F]$-module.
If $\chi$ denotes its character, we have by Corollary~\ref{IrredMod} that
\begin{align*}
\chi(\Lambda^{[3,k]}) &=
\frac{1}{|F_g|} \sum_{y \in F_{g}} z_{3,k}(g,y) \; \eta(y)\\
&= \frac{1}{3} \sum_{i=0}^{2}z_{3,k}(g,c^{3i}) \; \zeta^i = 1+\zeta^k
\end{align*}
In particular, we see that~$\nu_3(\chi)$ is not an integer, in fact not even a real number.

\subsection[Noncocommutative Sweedler powers]{} \label{SweedlPowNonCocom}
The second counterexample that we consider is suggested by  Proposition~\ref{CoprodSweedl}, where we saw that certain Sweedler powers behave like cocommutative elements with respect to the action of characters. However, we will see now that they are not really cocommutative elements. In fact, we will see that already the second Sweedler power of the integral is in general not cocommutative. 

By definition, we have $x \in G_{2,1}(g,y)$ if and only if
$x^2=y$ and $g(x^{-1}.g)=1$, which means that $x^{-1}.g=g^{-1}$. By Proposition~\ref{Coeff}, this implies $x^2.g=g$, which is also easy to see directly. If $|F|$ is odd, $x$ and~$x^2$ generate the same cyclic subgroup, and we get that $x.g=g$. In this case, we therefore have
$x \in G_{2,1}(g,y)$ if and only if $x^2=y$, $x.g=g$, and $g^2=1$. The formula in Proposition~\ref{SweedlPowIntExa} therefore yields
$$\Lambda^{[2]}= 
\frac{1}{|F|} 
\sum_{\substack{g \in G,\, x \in F\\g^2=1,\; x.g=g}} b_g \o x^2$$
and consequently we have
$$\Delta(\Lambda^{[2]})= 
\frac{1}{|F|} \sum_{\substack{g,h \in G,\, x \in F\\(gh)^2=1,\; x.(gh)=gh}} 
(b_g \o x^2) \o (b_h \o x^2)$$
Suppose now that $G=S_5$\index{$S_5$}, the symmetric group on five letters, and let $F$ be cyclic of order~$3$. Choose a generator~$c$ of~$F$\index{$c$} and define an action of~$F$ on~$G$  by the requirement that the action of~$c$ be the conjugation by the 3-cycle $\tau:=(1,2,3)$\index{$\tau$}:
$$c.\sigma = \tau \sigma \tau^{-1}$$
Consider the elements $g:=(3,4)$ and $h:=(3,4,5)$. We have
$gh=(4,5)$ and $hg=(3,5)$, which means that $c.(gh) = gh$, but 
$c.(hg) \neq hg$. Therefore, the coefficient
of $(b_g \o c^2) \o (b_h \o c^2)$ in the above expansion of $\Delta(\Lambda^{[2]})$ is $1/|F|$, whereas the coefficient of $(b_h \o c^2) \o (b_g \o c^2)$ is
zero. In particular, $\Lambda^{[2]}$ is not cocommutative.

\subsection[Noncentral Sweedler powers]{} \label{SweedlPowNonCent}
The third counterexample that we give is suggested by  Proposition~\ref{SweedlPowCent}, where we saw that certain Sweedler powers of the integral are central. As we said there, this is not the case for all Sweedler powers, and we will now confirm this assertion by explicitly constructing such a noncentral Sweedler power.

With the notation used there, we get as in the proof of Proposition~\ref{SweedlPowIntExa} that
\begin{align*}
\Lambda^{[4,2]} &=
\frac{1}{|F|} \sum_{x \in F} \sum_{\substack{g_1,g_2,g_3,g_4 \in G \\
g_1 g_2 g_3 g_4 = 1}} 
(b_{g_1} b_{x.g_3} b_{x^2.g_2} b_{x^3.g_4})\o x^4 \\
&=
\frac{1}{|F|} \sum_{\substack{x \in F, \, g \in G \\
g (x^{-2}.g) (x^{-1}.g) (x^{-3}.g) = 1}} 
b_g \o x^4
\end{align*}
Suppose now that $G=S_8$\index{$S_8$}, the symmetric group on eight letters, and let $F$ be cyclic of order~$8$. Choose a generator~$c$ of~$F$\index{$c$} and define an action of~$F$ on~$G$  by the requirement that the action of~$c$ be the conjugation by the 8-cycle $\tau:=(1,2,3,4,5,6,7,8)$\index{$\tau$}, so that
$c.\sigma = \tau \sigma \tau^{-1}$.
Now, if $g$ is the 8-cycle $(1,3,5,8,2,4,6,7)$, we have
\begin{align*}
g & (c^{-2}.g) (c^{-1}.g) (c^{-3}.g) =\\
&(1,3,5,8,2,4,6,7)(7,1,3,6,8,2,4,5)(8,2,4,7,1,3,5,6)(6,8,2,5,7,1,3,4)=\id
\end{align*}
Furthermore, we have
$c^4.g = (5,7,1,4,6,8,2,3) \neq g$.

Suppose now that $\Lambda^{[4,2]}$ were central. Then we would have
$(b_g \o 1)\Lambda^{[4,2]} 
= (b_g \o 1) \Lambda^{[4,2]} (b_g \o 1)$.
But from the formula above we get
\begin{align*}
&(b_g \o 1)\Lambda^{[4,2]} = 
\frac{1}{8} \sum^7_{\substack{i=0 \\
g (c^{-2i}.g) (c^{-i}.g) (c^{-3i}.g) = \id}} 
b_g \o c^{4i} \\
&(b_g \o 1) \Lambda^{[4,2]} (b_g \o 1)
= 
\frac{1}{8} \sum^7_{\substack{i=0 \\
g (c^{-2i}.g) (c^{-i}.g) (c^{-3i}.g) = \id\\
c^{4i}.g=g}} 
b_g \o c^{4i}
\end{align*}
Now the preceding calculations show that the coefficient of~$b_g \o c^4$
in the second sum is zero, whereas in the first sum the summand with~$i=1$
contributes to this coefficient, which implies, as we are in characteristic zero, that this coefficient is nonzero. We have therefore reached a contradiction, establishing that $\Lambda^{[4,2]}$ is not central.


\end{document}